\documentclass{article}
\usepackage{latexsym,amsmath,amssymb,amscd}
\begin{document}

\title{On Spherically Symmetric Motions of a Gaseous Star Governed by the Euler-Poisson Equations}
\author{Tetu Makino\footnote{Faculty of Engineering,
Yamaguchi University, Ube, 755-8611, Japan. E-mail : makino@yamaguchi-u.ac.jp}}
\date{\today}
\maketitle

\begin{abstract}
We consider spherically symmetric motions of a 
polytropic gas under the self-gravitation governed by
the Euler-Poisson equations. The adiabatic exponent (=
the ratio of the specific heats ) $\gamma$ is assumed to satisfy
$6/5< \gamma \leq 2$. Then there are equilibria touching the vacuum with finite radii, and 
the linearized equation around one of the equilibria 
has time-periodic solutions. To justify the linearization, we should construct true solutions
for which this time-periodic solution plus the equilibrium is the first approximation. 
We solve this problem
by the Nash-Moser theorem. The result will realize the so-called
physical vacuum boundary. But the present study restricts $\gamma$ to
the case in which $\gamma/(\gamma-1)$ is an integer. Other cases are reserved to 
the future as an open problem. The time-local existence
of smooth solutions to the Cauchy problems is also dicussed.


\end{abstract}

\newtheorem{Lemma}{Lemma}
\newtheorem{Proposition}{Proposition}
\section{Introduction}

We consider spherically symmetric motions of a gaseous star governed
by the Euler-Poisson equations:

\begin{eqnarray}
&& \frac{\partial\rho}{\partial t}+u\frac{\partial\rho}{\partial r}+\rho \frac{\partial u}{\partial r}+\frac{2}{r}\rho u=0, \nonumber \\
&&\rho\Big(\frac{\partial u}{\partial t}+u\frac{\partial u}{\partial r}\Big)
+\frac{\partial P}{\partial r} =-\rho\frac{\partial\Phi}{\partial r} \qquad (0<t, 0<r), \nonumber \\
&&\frac{1}{r^2}\frac{\partial}{\partial r}\Big(r^2\frac{\partial\Phi}{\partial r}\Big)=4\pi g_0 \rho.
\end{eqnarray}
Here $\rho$ is the density, $u$ the velocity, $P$ the pressure, $\Phi$ the gravitational potential,
and $g_0$ is the gravitational constant. In this work
we assume
\begin{equation}
P=A\rho^{\gamma},
\end{equation}
where $A$ and $\gamma$ are positive constants, and we assume $1<\gamma\leq 2$. \\

Introducing the mass
$$m:=4\pi\int_0^r\rho(t,r')r'^{2}dr', $$
we can write the equations as
\begin{eqnarray}
&&\frac{\partial\rho}{\partial t}+u\frac{\partial\rho}{\partial r}+\rho \frac{\partial u}{\partial r}+\frac{2}{r}\rho u=0, \nonumber \\
&&\rho\Big(\frac{\partial u}{\partial t}+u\frac{\partial u}{\partial r}\Big)
+\frac{\partial P}{\partial r}=-g_0\frac{\rho m}{r^2}.
\end{eqnarray}
\\

On the other hand, equilibria for the equations (1) are governed by the ordinary differential equation
$$-\frac{1}{r^2}\frac{d}{dr}\Big(\frac{r^2}{\rho}\frac{dP}{dr}\Big)=4\pi g_0\rho.$$

In order to normalize this equation, we put
$$\rho=\rho_c\theta^{\frac{1}{\gamma-1}}$$
and
$$r=\rho_c^{\frac{\gamma-2}{2}}K^{-1/2}\xi 
\qquad \mbox{with}\quad K:=\frac{4\pi g_0(\gamma-1)}{A\gamma},$$
where $\rho_c$ is an arbitrary positive number, say, the central density.
Then the equation for equilibria turns out to be
$$\frac{1}{\xi^2}\frac{d}{d\xi}\xi^2\frac{d\theta}{d\xi}+\theta^{\frac{1}{\gamma-1}}=0,
$$
which is called the `Lane-Emden equation'. The solution
$\theta(\xi)$ of the equation such that
$$\theta|_{\xi=0}=1,\qquad
\frac{d\theta}{d\xi}\Big|_{\xi=0}=0$$
is called the `Lane-Emden function of polytropic index $\frac{1}{\gamma-1}$'.
It is known that if and only if $6/5<\gamma$ there is a finite $\xi_1$
such that $\theta(\xi)>0$ for $0\leq \xi<\xi_1$ and
$\theta(\xi_1)=0$, and the radius $R$ and the total mass 
$$M:=4\pi\int_0^R\rho(r)r^2dr$$
of the equilibrium $\rho(r)$
are given by
$$R=\rho_c^{\frac{\gamma-2}{2}}K^{-1/2}\xi_1, 
\quad \mbox{and}\quad
M=4\pi \rho_c^{\frac{3\gamma-4}{2}}K^{-3/2}
\Big(-\xi^2\frac{d\theta}{d\xi}\Big)_{\xi=\xi_1}.$$
A numerical table of $\xi_1, (-\xi^2d\theta/d\xi)_{\xi=\xi_1}$ 
for various $\gamma$ can be found in \cite[p.96]{Chandrasekhar}. 

Anyway we have

\begin{Lemma}
Assume $6/5 < \gamma \leq 2$. For any positive number $\rho_c$ given, there is an equilibrium
$\rho=\bar{\rho}(r)$ with positive numbers $R, \rho_1$ such that
$\bar{\rho}(r)$ is positive and analytic in $0<r<R$ and
\begin{eqnarray*}
\bar{\rho}(r)&=& \rho_c(1+[r^2]_1) \qquad \mbox{as}\ r\rightarrow 0, \\
\bar{\rho}(r)&=&\rho_1(R-r)^{\frac{1}{\gamma-1}}(1+[R-r, (R-r)^{\frac{\gamma}{\gamma-1}}]_1) \qquad
\mbox{as}\ r\rightarrow R-0.
\end{eqnarray*}  
\end{Lemma}

{\bf Notational Remark} 
Here and hereafter $[X]_q$ denotes a power series of the form
$\sum_{j\geq q}a_jX^j$ with positive radius of convergence, and
$[X,Y]_q$ a convergent power series of the form
$\sum_{j+k\geq q}a_{jk}X^jY^k$. \\

For a proof of Lemma 1, see, e.g., \cite{JosephLundgren}, and \cite[Chapter V]{Lefschetz}
 or \cite[Chapter IX]{Wasow} and Appendix 1.\\ 

{\bf Remark} In the expansion of
$\bar{\rho}(r)$ as $r \rightarrow R$, the terms including \\
$(R-r)^{\frac{\gamma}{\gamma-1}}$ actually appear if $\displaystyle \frac{\gamma}{\gamma-1}$
is not an integer. Let us prove it.
Otherwise we would have
$$\bar{\rho}(r)=\rho_1(R-r)^{\frac{1}{\gamma-1}}(1+[R-r]_1) $$
and the function
$$U(r):=\bar{\rho}(r)^{\gamma-1}=\rho_1^{\gamma-1}(R-r)(1+[R-r]_1)$$
would be analytic at $r=R$. Now $U$ satisfies
$$-\frac{d^2U}{dr^2}-\frac{2}{r}\frac{dU}{dr}=KU^{\frac{1}{\gamma-1}},
\quad K:=\frac{4\pi g_0(\gamma-1)}{A\gamma}.$$
Since $U$ is analytic, the left-hand side is analytic, and so, the right-hand side 
$$K\rho_1 (R-r)^{\frac{1}{\gamma-1}}(1+[R-r]_1)$$
would be analytic at $r=R$. Then $\displaystyle \frac{1}{\gamma-1}$ shoud be an integer. This contradicts
to that $\displaystyle \frac{\gamma}{\gamma-1}=\frac{1}{\gamma-1}+1$ is not an integer. 

In fact we can find that, if $\displaystyle \frac{\gamma}{\gamma-1}\not\in\mathbb{N}$, then
\begin{eqnarray*}
\bar{\rho}^{\gamma-1}&=&U=
C(R-r)\Big(1+\\
&&+\frac{1}{R}(R-r)-
\frac{(\gamma-1)^2KC^{\frac{2-\gamma}{\gamma-1}}}{\gamma(2\gamma-1)}(R-r)^{\frac{\gamma}{\gamma-1}}+
[R-r, (R-r)^{\frac{\gamma}{\gamma-1}}]_2\Big) 
\end{eqnarray*}
and
\begin{eqnarray*}\frac{1}{\bar{\rho}}\frac{d\bar{P}}{dr}=\frac{A\gamma}{\gamma-1}\frac{dU}{dr}&=&
-\frac{A\gamma}{\gamma-1}C\Big(
1+\frac{2}{R}(R-r)+ \\
&&-\frac{(\gamma-1)KC^{\frac{2-\gamma}{\gamma-1}}}{\gamma}(R-r)^{\frac{\gamma}{\gamma-1}}+
[R-r, (R-r)^{\frac{\gamma}{\gamma-1}}]_2\Big),
\end{eqnarray*}
where
$C=\rho_1^{\gamma-1}$
and $\bar{P}(r)=A\bar{\rho}(r)^{\gamma}$.\\

%\subsection{}
\textbullet\  In the following discussion we assume that $6/5<\gamma\leq 2$ and we fix an equilibrium $\bar{\rho}(r)$
with the properties in the above Lemma.\\

We are going to construct solutions around this fixed equilibrium. \\

Here let us glance at the history of researches of this problem.

Of course there were a lot of works on the Cauchy problem to the
compressible Euler equations. But there were gaps if we consider 
density distributions which contain vacuum regions.

As for local-in-time existence of smooth density
with compact support, \cite{M1986} treated the problem under the assumption
that the initial density is non-negative and the initial value of
$$\omega:=\frac{2\sqrt{A\gamma}}{\gamma-1}\rho^{\frac{\gamma-1}{2}}$$
is smooth, too. By the variables $(\omega,u)$ the equations are symmetrizable
continuously including the region of vacuum. Hence the theory of 
quasi-linear symmetric hyperbolic
systems can be applied. 
However, since
$$\omega\propto \Big(\frac{1}{r}-\frac{1}{R}\Big)^{\frac{1}{2}}\sim 
\mbox{Const.}(R-r)^{\frac{1}{2}} \quad \mbox{as}\  r\rightarrow R-0$$
for equilibria, $\omega$ is not smooth at the boundary $r=R$ with the vacuum.
Hence the class of solutions considered in \cite{M1986} cannot cover equilibria.
( See \cite{MakinoUkai} for the discussion on non-isentropic cases. The situation is similar.)

On the other hand, possibly discontinuous weak solutions with compactly
supported density can be constructed. The article \cite{M1997} gave local-in-time
existence of bounded weak solutions under the assumption
that the initial density is bounded and non-negative, provided that
the gas is confined to the domain outside a solid ball. The proof by the 
compensated compactness method is due to \cite{MT}, and 
\cite{DCL}. Of course the class of weak solutions can cover equilibria, 
but the concrete structures of solutions were not so clear. 

Therefore we wish to construct solutions whose regularities 
are weaker than solutions with smooth $\omega$ and stronger than possibly discontinuous 
weak solutions. The present result is an answer to this wish. 
More concretely speaking, the solution $(\rho(t,r),u(t,r))$ constructed
in this article should be continuous on
$0\leq t\leq T,0\leq r <\infty$ and there should be found a continuous curve
$r=R_F(t), 0\leq t\leq T,$ such that
$|R_F(t)-R|\ll 1, \rho(t,r)>0 $ for $  0\leq t\leq T, 0\leq r <R_F(t)$
and $\rho(t,r)=0$ for $0\leq t\leq T, R_F(t)\leq r<\infty$. The curve $r=R_F(t)$ is
the free boundary at which the density touches the vacuum.
It will be shown that the solution satisfies
$$\rho(t,r)=C(t)(R_F(t)-r)^{\frac{1}{\gamma-1}}(1+O(R_F(t)-r))$$
as $r \rightarrow R_F(t)-0$. Here $C(t)$ is positive and smooth in $t$.
This situation is ``physical vacuum boundary" so-called by
\cite{JM} and \cite{CS}. This concept can be traced back to
\cite{L}, \cite{LY}, \cite{Y}. Of course this singularity is just that
of equilibria.\\

Since the major difficulty
comes from the free boundary touching the vacuum, which moves along time. So,
we take the Lagrangian mass coordinate $m$ as the independent variable
instead of $r$.
Then we can write the equations as
\begin{eqnarray*}
&&\frac{\partial\rho}{\partial t}+4\pi\rho^2(r^2u)_m=0, \\
&&\frac{\partial u}{\partial t}+4\pi r^2P_m=-g_0\frac{m}{r^2}, \\
&&r=\Big(\frac{3}{4\pi}\int_0^m\frac{dm}{\rho}\Big)^{1/3}.
\end{eqnarray*}

Since
$$\frac{\partial r}{\partial t}=u, \qquad \frac{\partial r}{\partial m}=\frac{1}{4\pi\rho r^2}, $$
the equations are reduced to the single second order equation
\begin{equation}
r_{tt}+4\pi r^2P_m=-g_0\frac{m}{r^2},
\end{equation}
where
$$P=A\Big(4\pi r^2\frac{\partial r}{\partial m}\Big)^{-\gamma}.$$\\

Now we derive the equation for the perturbation $y$ defined by
\begin{equation}
r(t,m)=\bar{r}(m)(1+y(t,\bar{r}(m))).
\end{equation}
Here $m\mapsto \bar{r}(m)$ is the function of the Lagrangian mass variable $m$
associated with the fixed equilibrium. In other words,
it is the inverse function of
$$\bar{r} \mapsto m=4\pi\int_0^{\bar{r}}\bar{\rho}(r')r'^2dr'.$$

Keeping in mind
$$\frac{\partial r}{\partial m }=
\frac{\partial \bar{r}}{\partial m}\Big(1+y+\frac{\bar{r}}{\bar{r}_m}\frac{\partial y}{\partial m}\Big),$$
we have
$$P=\bar{P}\Big(1-G(y, \frac{\bar{r}}{\bar{r}_m}\frac{\partial y}{\partial m})
\Big). $$
Here $G(y,v)=3\gamma y+\gamma v+[y,v]_2$ is defined by
$$(1+y)^{-2\gamma}(1+y+v)^{-\gamma}=1-G(y,v).$$

Then the equation is reduced to
$$\bar{r}y_{tt}
+\frac{1}{\bar{\rho}}(1+y)^2
\frac{\partial}{\partial \bar{r}}\Big(\bar{P}\Big(
1-G(y,\bar{r}\frac{\partial y}{\partial \bar{r}})\Big)\Big) +
g_0\frac{m}{\bar{r}^2(1+y)^2}=0, $$
where we have used
$$\frac{\partial}{\partial m}=\bar{r}_m\frac{\partial}{\partial \bar{r}}=
\frac{1}{4\pi\bar{\rho}\bar{r}^2}\frac{\partial}{\partial\bar{r}}. $$
We note that the equilibrium satisfies 
$$\frac{1}{\bar{\rho}}\frac{\partial\bar{P}}{\partial\bar{r}}+g_0\frac{m}{\bar{r}^2}=0.$$
Let us introduce $H(y)=4y+[y]_2$ by
$$H(y)=(1+y)^2-\frac{1}{(1+y)^2}.$$
Then the equation can be written as
\begin{equation}
\frac{\partial^2y}{\partial t^2}
-\frac{1}{\rho r}(1+y)^2\frac{\partial}{\partial r}\Big(PG(y,r\frac{\partial y}{\partial r})\Big)
+\frac{1}{\rho r}\frac{dP}{dr}H(y)=0.
\end{equation}

Here we have used the abbreviations 
$r, \rho, P, \frac{dP}{dr}$
instead of
$\bar{r}, \bar{\rho}, \bar{P}, \frac{d\bar{P}}{d\bar{r}}$. We consider this nonlinear wave equation.

It is easy to verify by a scale transformation of variables
that we can assume that $A=1/\gamma$ so that $P=\rho^{\gamma}/\gamma $ 
without loss of generality. Hence we assume so. \\

Here let us propose the main goal of this study roughly.
Let us fix an arbitrarily large positive number $T$. Then,
under the condition that $\gamma/(\gamma-1)$ is an integer,  we have \\

{\bf Main Goal } {\it For sufficiently small $\varepsilon>0$ there
is a solution $y=y(t,r;\varepsilon)$ of (6) in
$C^2([0,T]\times[0,R])$ such that
$$y(t,r;\varepsilon)=\varepsilon y_1(t,r)+O(\varepsilon^2).$$
The same estimates $O(\varepsilon^2)$ hold between the
higher order derivatives of $y$ and $\varepsilon y_1$.}\\

Here $y_1(t,r)$ is a time-periodic function specified in Section 2, which
is of the form
$$y_1(t,r)=\sin(\sqrt{\lambda}t+\theta_0)\cdot \Phi(r),$$
where $\lambda$ is a positive number, $\theta_0$ a constant, and
$\Phi(r)$ is an analytic function of $0\leq r\leq R$.\\

Once the solution $y(t,r;\varepsilon)$ is given, then
 the corresponding motion of gas particles can be expressed by the Lagrangian 
coordinate as
\begin{eqnarray*}
r(t,m)&=&\bar{r}(m)(1+y(t,\bar{r}(m);\varepsilon)) \\
&=&\bar{r}(m)(1+\varepsilon y_1(t,\bar{r}(m))+O(\varepsilon^2)).
\end{eqnarray*}
The curve $r=R_F(t)$ of the free vacuum boundary is given by
$$R_F(t)=r(t,M)=R(1+\varepsilon \sin(\sqrt{\lambda}t+\theta_0)\Phi(R)+O(\varepsilon^2)).$$
{\it The free boundary $R_F(t)$ oscillates around $R$ with time-period $2\pi/\sqrt{\lambda}$ approximately.}

The solution $(\rho,u)$ of the original problem (1)(2)
is given by
$$\rho=\bar{\rho}(\bar{r})\Big((1+y)^2\Big(1+y+\bar{r}\frac{\partial y}{\partial\bar{r}}\Big)\Big)^{-1},
\qquad u=\bar{r}\frac{\partial y}{\partial t} $$
implicitly by
\begin{eqnarray*}
\bar{r}&=&\bar{r}(m), \qquad y=y(t,\bar{r}(m);\varepsilon) \\
\frac{\partial y}{\partial\bar{r}}&=&
\partial_ry(t,\bar{r}(m);\varepsilon), \qquad
\frac{\partial y}{\partial t}=
\partial_ty(t,\bar{r}(m);\varepsilon),
\end{eqnarray*}
where $m=m(t,r)$ for $0\leq r\leq R_F(t)$. Here
$r\mapsto m=m(t,r)$ is given as the inverse function of the 
function
$$ m\mapsto r=r(t,m)=\bar{r}(m)(1+y(t,\bar{r}(m);\varepsilon)).$$
We note that
$$R_F(t)-r(t,m)=R(1+y(t,R;\varepsilon))-\bar{r}(m)(1+y(t,\bar{r}(m)
;\varepsilon))$$
implies
$$\frac{1}{\kappa}(R-\bar{r})\leq R_F(t)-r\leq \kappa (R-\bar{r})
$$
with $0<\kappa-1\ll 1$, since
$|y|+|\partial_r y|\leq \varepsilon C$. Therefore
$$y(t,\bar{r}(m);\varepsilon)=y(t,R; \varepsilon)+O(R_F(t)-r(t,m)),$$
and so on.
Hence we get the ``physical vacuum boundary", that is,
the corresponding density distribution $\rho=\rho(t,r)$,
where $r$ denotes the original Euler coordinate,
satisfies
$$\rho(t,r)>0\  \mbox{for}\  0\leq r<R_F(t), \qquad
\rho(t,r)=0\  \mbox{for}\  R_F(t)\leq r, $$
and,
since $y(t,r)$ is smooth on $0\leq r\leq R$, we have
$$\rho(t,r)=C(t)(R_F(t)-r)^{\frac{1}{\gamma-1}}(1+O(R_F(t)-r)) $$
as $r \rightarrow R_F(t)-0$. Here $C(t)$ is positive and smooth in $t$.

\section{Analysis of the linearized equation}

The linearized equation is
\begin{eqnarray}
&&\frac{\partial^2y}{\partial t^2}+\mathcal{L}y=0, \\
\mathcal{L}y&:=&
-\frac{1}{\rho r}\frac{\partial }{\partial r}
\Big(P\Big(3\gamma y+\gamma r\frac{\partial y}{\partial r}\Big)\Big)+
\frac{1}{\rho r}\frac{dP}{dr}\cdot(4y) \nonumber \\
&=&-\frac{1}{\rho r^4}\frac{\partial}{\partial r}\Big(\gamma r^4P\frac{\partial y}{\partial r}\Big)
+\frac{1}{\rho r}(4-3\gamma)\frac{dP}{dr}y,
\end{eqnarray}
and the associated eigenvalue problem is
$\mathcal{L}y=\lambda y.$\\

This eigenvalue problem was first wrote down in
\cite[p.10, (12)]{Eddington}(1918).
But the spectral property of the operator, whose coefficients are singular,
had been long believed as a Sturm-Liouville type without proof. A mathematically rigorous discussion
was first done by \cite{Beyer}(1995). The essential point is as follows.\\

Let us use the Liouville transformation:

$$\xi:=\int_0^r\sqrt{\frac{\rho}{\gamma P}}dr, \qquad \eta:= r^2(\gamma P\rho)^{\frac{1}{4}}y. $$

Through this transformation the equation
$$\mathcal{L}y=\lambda y+f$$
turns out to be the standard form
$$-\frac{d^2\eta}{d\xi^2}+q\eta=\lambda\eta + \hat{f}, $$
where
\begin{eqnarray*}
q&=&\frac{\gamma P}{\rho}\Big(
\frac{2}{r^2}+\Big(\frac{7-3\gamma}{2}+\frac{1+\gamma}{4}
\frac{rm_r}{m}\Big)\frac{1}{r\rho}\frac{d\rho}{dr}+
\frac{(\gamma+1)(3-\gamma)}{16}\Big(\frac{1}{\rho}\frac{d\rho}{dr}\Big)^2\Big), \\
\hat{f}&=&r^2(\gamma P\rho)^{\frac{1}{4}}f.
\end{eqnarray*}
The variable $\xi$ runs on the interval $(0, \xi_+)$, where
$$\xi_+:=\int_0^R\sqrt{\frac{\rho}{\gamma P}}dr< \infty.$$

Since
$$\xi \sim \sqrt{\frac{\rho_c}{\gamma P_c}}r \qquad \mbox{as}\ r\rightarrow 0, $$
we see
$$q \sim \frac{\gamma P_c}{\rho_c}\frac{2}{r^2} \sim \frac{2}{\xi^2} \qquad \mbox{as}\ \xi\rightarrow 0.$$

Since
$$\frac{1}{\rho}\frac{d\rho}{dr}\sim-\frac{1}{\gamma-1}(R-r)^{-1}, \qquad \frac{\gamma P}{\rho}\sim
\rho_1^{\gamma-1}(R-r) \qquad \mbox{as}\ r\rightarrow R, $$
and
$$R-r \sim \frac{1}{4}\rho_1^{\gamma-1}(\xi_+-\xi)^2 \qquad\mbox{as}\  \xi \rightarrow \xi_+,$$
we see
$$q\sim \frac{\gamma P}{\rho}\frac{(\gamma+1)(3-\gamma)}{16}\Big(\frac{1}{\rho}\frac{d\rho}{dr}\Big)^2
\sim
\frac{1}{4}\frac{(1+\gamma)(3-\gamma)}{(\gamma-1)^2}\frac{1}{(\xi_+-\xi)^2}$$
as $\xi\rightarrow \xi_+$.  It follows from $1<\gamma<2$ that 
$$\frac{1}{4}\frac{(1+\gamma)(3-\gamma)}{(\gamma-1)^2}>\frac{3}{4}. $$

Of course $q$ is bounded from below, but it is difficult to know 
whether its minimum is positive or not. Anyway, the both boundary points 
$\xi=0, \xi_+$ are of limit point type, provided that $1<\gamma<2$. See, e.g., \cite[p.159, Theorem X.10]{Reed}.
The exceptional case $\gamma=2$ 
will be discussed later.
See the discussion after Lemma 2 below.
Hence we have the following conclusion:

\begin{Proposition}
The operator $\mathfrak{T}_0, \mathcal{D}(\mathfrak{T}_0)=
C_0^{\infty}(0, \xi_+), \mathfrak{T}_0\eta =-\eta_{\xi\xi}+q\eta, $
in $L^2(0,\xi_+)$ has the Friedrichs extension $\mathfrak{T}$, a self-adjoint operator, whose spectrum
consists of simple eigenvalues
$\lambda_1<\cdots<\lambda_n<\lambda_{n+1}<\cdots\rightarrow +\infty$. 
In other words,
the operator
$\mathfrak{S}_0, \mathcal{D}(\mathfrak{S}_0)=C_0^{\infty}(0,R), \mathfrak{S}_0y=\mathcal{L}y$ in
$L^2((0,R), r^4\rho dr)$ has the Friedrichs extension $\mathfrak{S}$, a self-adjoint 
operator with eigenvalues $(\lambda_n)_n$. 
\end{Proposition}

The domain $\mathcal{D}(\mathfrak{T})$ of the Firedrichs extension $\mathfrak{T}$ is, by
definition, 
\begin{eqnarray*}
\mathcal{D}(\mathfrak{T})=\{\eta \in L^2(0,\xi_+)&|& \exists \phi_n\in C_0^{\infty}(0,\xi_+), 
\quad Q[\phi_m-\phi_n]\rightarrow 0 \\
&&\mbox{as}\  m, n\rightarrow \infty, \quad \phi_n\rightarrow \eta \quad\mbox{in}\  L^2(0,\xi_+) \\
&&\mbox{and}\  -\eta_{\xi\xi}+q\eta \in L^2(0,\xi_+) \quad\mbox{in distribution sense}\},
\end{eqnarray*}
where
$$Q[\phi]:= \int_0^{\xi_+}\Big(\Big|\frac{d\phi}{d\xi}\Big|^2+(q+c)|\phi|^2\Big)d\xi,$$
and
$c$ is a constant $>|\min q|$. But $\mathcal{D}(\mathfrak{T})$ is characterized as follows:
$$\mathcal{D}(\mathfrak{T})=\{\eta \in C[0,\xi_+]\quad|\quad 
\eta(0)=\eta(\xi_+)=0,\quad -\eta_{\xi\xi}+q\eta \in L^2(0,\xi_+) \}.$$

Let us prove it, denoting by $M$ the right-hand side.
Let $\eta \in \mathcal{D}(\mathfrak{T})$. Then there are $\phi_n \in C_0^{\infty}(0,\xi_+)$ such
that $\phi_n \rightarrow \eta$ in $L^2$ and $Q[\phi_m-\phi_n]\rightarrow 0$.
Since
$$|\phi_m(\xi)-\phi_n(\xi)|\leq \sqrt{\xi}\Big(\int_0^{\xi}((\phi_m-\phi_n)_{\xi})^2d\xi\Big)^{1/2}
\leq \sqrt{\xi}(Q[\phi_m-\phi_n])^{1/2}\rightarrow 0,$$
we have $\phi_n \rightarrow \eta $ uniformly on $[0,\xi_+]$. Hence
$\eta\in C[0,\xi_+]$ and $\eta(0)=0$. Similarly $\eta(\xi_+)=0$. Thus 
$\mathcal{D}(\mathfrak{T}) \subset M$. Let $\eta \in M$. Put
$f:=-\eta_{\xi\xi}+q\eta \in L^2$. Then 
$-\eta_{\xi\xi}+(q+c)\eta = g:=f+c\eta \in L^2$.
Since $0$ belongs to the resolvent set of $\mathfrak{T}+c$, we have $v:=(\mathfrak{T}+c)^{-1}g
\in \mathcal{D}(\mathfrak{T})$. Hence $w:=\eta-v \in C[0,\xi_+]$ and
$w(0)=w(\xi_+)=0, -w_{\xi\xi}+(q+c)w=0$, for $\mathcal{D}(\mathfrak{T}) \subset M$.
Using $q+c >0$, we can deduce that $w \equiv 0$ and $\eta=v \in \mathcal{D}(\mathfrak{T})$,
that is, $M\subset \mathcal{D}(\mathfrak{T})$. (In fact,
if $w$ did not vanish identically, there would exist
$a \in (0,\xi_+)$ such that
$Dw(a)=0$ and $w(a)\not= 0$.
If $w(a)>0$, then 
$$Dw(\xi)=\int_a^{\xi} D^2w(\xi')d\xi' =
\int_a^{\xi}(q+c)w(\xi')d\xi'
$$
implies $Dw(\xi)>0$ for $a < \xi <\xi_+$ and it contradicts to
$w(\xi_+)=0$.
If $w(a)<0$, then 
$$Dw(\xi)=-\int_{\xi}^a(q+c)w(\xi')d\xi' $$
implies $Dw(\xi)>0$ for $0<\xi <a$ and it contradicts to $w(0)=0$.)
$\blacksquare$ \\

Although it is not easy to judge the signature of $\min q$, we have

\begin{Proposition}[\cite{SSLin},1997]
If and only if $ 4/3 < \gamma \leq 2$, the least eigenvalue $\lambda_1$ is positive.
\end{Proposition}

{\bf Proof}\  The function $y\equiv 1$ satisfies
$$\mathcal{L}y=\frac{1}{\rho r}(4-3\gamma)\frac{dP}{dr} =:f>0. $$
Let us consider the corresponding function
$$\eta_1=r^2(\gamma P\rho)^{\frac{1}{4}}$$ through the Liouville transformation.
It is easy to show that $\eta_1$ and $d\eta_1/d\xi$ vanish at $\xi=0, \xi_+$ and $\eta_1 \in \mathcal{D}(\mathfrak{T})$.
Let $\phi_1(\xi)$ be the eigenfunction of 
$-d^2/d\xi^2 +q$ associated with the least eigenvalue $\lambda_1$. We can assume
$\phi_1(\xi)>0$ for $0<\xi<\xi_+$ and 
$\phi_1$ and $d\phi_1/d\xi$ vanishes at $\xi=0, \xi_+$. Then the integration by parts gives
$$\lambda_1\int_0^{\xi_+}\phi_1\eta_1d\xi=
\int_0^{\xi_+}\phi_1(-\eta_{1,\xi\xi}+q\eta_1)d\xi . $$
Since
$$-\eta_{1,\xi\xi}+q\eta_1=\hat{f}=r(\gamma P)^{\frac{1}{4}}\rho^{-\frac{3}{4}}(4-3\gamma)\frac{dP}{dr}  $$
and $\frac{dP}{dr} <0$, we have the assertion. $\blacksquare$\\

{\bf Remark}\  Assume that $3/4 <\gamma \leq 2$. Then the least eigenvalue, 
which is positive, is given by the variational formula
$$\lambda_1=\min \frac{(\mathcal{L}y|y)_{\mathfrak{X}}}{\|y\|_{\mathfrak{X}}^2},$$
where $\mathfrak{X}=L^2((0,R),\rho r^4dr))$ endowed with
 $(u|v)_{\mathfrak{X}}=\int_0^R uv\rho r^4dr$. From this
we can deduce the following
Ritter-Eddington's law of the period-density relation: {\it Let us
consider
equilibria $\rho(r)$ with $\rho(0)=\rho_c$ and the corresponding least eigenvalue $\lambda_1$ or
the ``period" $\Pi:=2\pi/\sqrt{\lambda_1}$; then $\Pi\sqrt{\rho_c}$ is a constant depending 
only upon $g_0, A, \gamma$.}

In fact we can consider the one parameter
family of equilibria
$$\rho(r)=\rho_\kappa(r):=\kappa^{\frac{2}{\gamma-2}}\bar{\rho}(r/\kappa)$$
which has radius $R=\kappa \bar{R}\propto \kappa$ and the central density
$\rho_c=\kappa^{\frac{2}{\gamma-2}}\bar{\rho}_c\propto \kappa^{\frac{2}{\gamma-2}}$. Here $\bar{\rho}$ is a fixed equilibrium with radius $\bar{R}$
and central density $\bar{\rho}_c$. Then it is easy to see
that
$(\mathcal{L}y|y)_{\mathfrak{X}}= \kappa^{\frac{5\gamma-6}{\gamma-2}}(\bar{\mathcal{L}}y_{\kappa}|y_{\kappa})_{\bar{\mathfrak{X}}}$,
where $y_{\kappa}(\bar{r})=y(\kappa\bar{r})$ and
$\bar{\mathfrak{X}}=L^2((0,\bar{R}),\bar{\rho}\bar{r}^4d\bar{r})$, and
$\|y\|_{\mathfrak{X}}^2= \kappa^{\frac{5\gamma-8}{\gamma-2}}\|y_{\kappa}\|_{\bar{\mathfrak{X}}}^2$. Hence
we have $\lambda_1\propto \kappa^{\frac{2}{\gamma-2}}\propto \rho_c$. This
completes the proof. 
( Note that the mean density $M/(4\pi R^3/3) \propto \kappa^{\frac{2}{\gamma-2}}\propto\   
\mbox{the central density}\ \rho_c$.)
 This fact was stated in \cite[p. 15]{Eddington},
 as a result that the pulsation theory conforms with observation of variable stars. 
As for the priority of A. Ritter (1879), see \cite{Rosseland}.\\

Let us introduce the variable $x$ defined by
\begin{equation}
x:=\frac{\tan^2\theta}{1+\tan^2\theta}, \qquad \theta:=\frac{\kappa\xi}{2}=
\frac{\kappa}{2}\int_0^r\sqrt{\frac{\rho}{\gamma P}}dr,
\end{equation}
where $\kappa=\pi/\xi_+$. Then $x$ runs over the interval $[0,1]$ while $r$
runs over $[0,1]$, and
$$\frac{dx}{dr}=
\kappa\sqrt{x(1-x)}\sqrt{\frac{\rho}{\gamma P}}=
\kappa\sqrt{x(1-x)}\rho^{\frac{-\gamma+1}{2}}.$$
Since
\begin{eqnarray*}
\frac{d}{dr}&=&\kappa\sqrt{x(1-x)}\rho^{\frac{-\gamma+1}{2}}\frac{d}{dx}, \\
\frac{d^2}{dr^2}&=&\kappa^2x(1-x)\rho^{-\gamma+1}\frac{d^2}{dx^2}+ \\
&&+\Big(\frac{1}{2}\kappa^2(1-2x)\rho^{-\gamma+1}+
\frac{-\gamma+1}{2}
\kappa\sqrt{x(1-x)}\rho^{\frac{-\gamma-1}{2}}\frac{d\rho}{dr}\Big)\frac{d}{dx},
\end{eqnarray*}
we have
\begin{eqnarray*}
\kappa^{-2}\mathcal{L}y&=&-x(1-x)\frac{d^2y}{dx^2}+ \\
&&-\Big(\frac{1}{2}(1-2x)+
\frac{\gamma+1}{2}\frac{1}{\kappa}\sqrt{x(1-x)}\rho^{\frac{\gamma-3}{2}}\frac{d\rho}{dr}+
\frac{4}{r}\frac{1}{\kappa}\sqrt{x(1-x)}\rho^{\frac{\gamma-1}{2}}\Big)\frac{dy}{dx}+ \\
&&+\frac{1}{\kappa^2}\frac{\rho^{\gamma-2}}{r}\frac{d\rho}{dr}(4-3\gamma)y.
\end{eqnarray*}\\

As $r \rightarrow 0 (x \rightarrow 0)$ we have
\begin{eqnarray*}
x&=&\frac{\kappa^2}{4}\rho_c^{-\gamma+1}r^2(1+[r^2]_1), \\
r&=&\frac{2}{\kappa}\rho_c^{\frac{\gamma-1}{2}}\sqrt{x}(1+[x]_1), \\
\frac{d\rho}{dr}&=&r[r^2]_0, \\
\frac{4}{r}\frac{1}{\kappa}\sqrt{x(1-x)}\rho^{\frac{\gamma-1}{2}}&=&2+[x]_1.
\end{eqnarray*}
Then it follows that
$$\kappa^{-2}\mathcal{L}y=
-x(1-x)\frac{d^2y}{dx^2}-\Big(
\frac{5}{2}+[x]_1\Big)\frac{dy}{dx}+[x]_0y.
$$\\

On the other hand, as $r \rightarrow R (x \rightarrow 1)$, we have
\begin{eqnarray*}
1-x&=&\kappa^2\rho_1^{-\gamma+1}(R-r)(1+[R-r,(R-r)^{\frac{\gamma}{\gamma-1}}]_1), \\
R-r&=&\frac{1}{\kappa^2}\rho_1^{\gamma-1}(1-x)(1+[1-x,(1-x)^{\frac{\gamma}{\gamma-1}}]_1), \\
\frac{d\rho}{dr}&=&-\frac{\rho_1}{\gamma-1}(R-r)^{\frac{2-\gamma}{\gamma-1}}(1+
[R-r,(R-r)^{\frac{\gamma}{\gamma-1}}]_1).
\end{eqnarray*}
Then it follows that
\begin{eqnarray*}
\kappa^{-2}\mathcal{L}y&=&
-x(1-x)\frac{d^2y}{dx^2} + \\
&&+\Big(\frac{\gamma}{\gamma-1}+[1-x,(1-x)^{\frac{\gamma}{\gamma-1}}]_1\Big)\frac{dy}{dx}
+[1-x,(1-x)^{\frac{\gamma}{\gamma-1}}]_0y.
\end{eqnarray*}

Changing the scale of $t$, we can and shall assume that $\kappa=1$
without loss of generality.\\

Summing up, we have:

\begin{Proposition}
We can write
\begin{equation}
\mathcal{L}y=-x(1-x)\frac{d^2y}{dx^2}-
\Big(\frac{5}{2}(1-x)-\frac{N}{2}x\Big)\frac{dy}{dx}+L_1(x)\frac{dy}{dx}+L_0(x)y,
\end{equation}
where
\begin{eqnarray*}
L_1(x)&=& 
\begin{cases}
[x]_1 & \quad\mbox{as}\quad x\rightarrow +0 \\
[1-x,(1-x)^{\frac{N}{2}}]_1& \quad\mbox{as}\quad x\rightarrow 1-0,
\end{cases}\\
L_0(x)&=&
\begin{cases}
[x]_0 &\quad \mbox{as}\quad x\rightarrow +0 \\
[1-x,(1-x)^{\frac{N}{2}}]_0 &\quad\mbox{as}\quad x\rightarrow 1-0
\end{cases}
\end{eqnarray*}
\end{Proposition}

Here $N$ is the parameter defined by
\begin{equation}
N=\frac{2\gamma}{\gamma-1} \quad
\Leftrightarrow \quad
\gamma=1+\frac{2}{N-2}.
\end{equation}
\\

Now let us fix a positive eigenvalue $\lambda=\lambda_n$ and an associated 
eigenfunction $\Phi(r)$ of $\mathcal{L}$. Then
$$y_1(t,r)=\sin(\sqrt{\lambda}t+\theta_0)\Phi(r)$$
is a time-periodic solution of the linearized problem.

Moreover we can claim 
\begin{Proposition}
We have
\begin{eqnarray*}
\Phi(r)&=&C_0(1+[r^2]_1) \qquad \mbox{as}\ r\rightarrow 0, \\
&=&C_0(1+[x]_1) \qquad \mbox{as}\ x\rightarrow 0
\end{eqnarray*}
and
\begin{eqnarray*}
\Phi(r)&=&C_1(1+[R-r, (R-r)^{\frac{\gamma}{\gamma-1}}]_1) \qquad \mbox{as}\ r\rightarrow R, \\
&=&C_1(1+[1-x, (1-x)^{\frac{N}{2}}]_1) \qquad \mbox{as}\ x\rightarrow 1
\end{eqnarray*}
Here $C_0$ and $C_1$ are non-zero constants. Other independent solutions of
$\mathcal{L}y=\lambda y$ do not belong to $L^2(r^4\rho dr)$ at $ r\sim R$. 
\end{Proposition}

To prove this,
we use the following lemma:

\begin{Lemma}
Let us consider the equation
$$z\frac{d^2y}{dz^2}+b(z,z^a)\frac{dy}{dz}=c(z,z^a)y, $$
where
$$b(z,z^{a})=a+[z,z^{a}]_1, \qquad c(z,z^{a})=[z,z^{a}]_0, $$
and let the positive number $a$ satisfy $a\geq 2$. Then 1) there is a solution $y_1$ of the form
$$y_1=1+[z,z^a]_1, $$
and 2) there is a solution $y_2$ of the form
$$y_2=z^{-a+1}(1+[z,z^a]_1) $$
provided $a \not\in \mathbb{N}$, or
$$y_2=z^{-a+1}(1+
[z, z^a]_1)+hy_1\log z $$
provided $a \in \mathbb{N}$. Here $h$ is a constant which can vanish in some
cases.
\end{Lemma}
For a proof, see \cite[Chapter 4]{CoddingtonLevinson}.

We apply this lemma for $a=\gamma/(\gamma-1)=N/2(\geq 2)$ and $z=1-x$.
Even if
$N=4 (\gamma=2)$, $y_2\sim
z^{-\frac{N}{2}+1}$ does not belong to $L^2(r^4\rho dr)=L^2(x^{3/2}(1-x)^{N/2-1}dx)$, and the boundary
point $r=R$ is of the limit point type.

\section{Statement of the main result}

We rewrite the equation (6) by using the linearized operator
$\mathcal{L}$ defined by (8) as
\begin{equation}
\frac{\partial^2y}{\partial t^2}+\Big(
1+G_I\Big(y, r\frac{\partial y}{\partial r}\Big)\Big)\mathcal{L}y+
G_{II}\Big(r,y,r\frac{\partial y}{\partial r}\Big)=0,
\end{equation}
where 
\begin{eqnarray*}
G_I(y,v)&=&(1+y)^2\Big(1+\frac{1}{\gamma}\partial_vG_2(y,v)\Big)-1, \\
G_{II}(r,y,v)&=&\frac{P}{\rho r^2}G_{II0}(y,v)
+\frac{1}{\rho r}\frac{dP}{dr}G_{II1}(y,v), \\
G_{II0}(y,v)&=&(1+y)^2(3\partial_vG_2-\partial_yG_2)v \\
&=&-2\gamma(1+y)^{-2\gamma+1}(1+y+v)^{-\gamma-1}v^2, \\
G_{II1}(y,v)&=&\frac{(1+y)^2}{\gamma}\partial_vG_2\cdot((-4+3\gamma)y+\gamma v) + \\
&&+H-4y(1+y)^2-(1+y)^2G_2 .
\end{eqnarray*}
Here
\begin{eqnarray*}
G_2(y,v)&:=&G(y,v)-(3\gamma y+\gamma v)=[y,v]_2, \\
\partial_vG_2&:=&\frac{\partial}{\partial v}G_2=\frac{\partial G}{\partial v}-\gamma
=[y,v]_1.
\end{eqnarray*}

We have fixed a solution $y_1$ of the linearized equation $y_{tt}+\mathcal{L}y=0$, and
we seek a solution $y$ of (6) or (12) of the form
$$y=\varepsilon y_1+\varepsilon w,$$
where $\varepsilon$ is a small positive parameter.
Then the equation which $w$ should satisfy turns out to be
\begin{eqnarray}
&&\frac{\partial^2w}{\partial t^2}+
\Big(1+\varepsilon a\Big(t,r,w,r\frac{\partial w}{\partial r},\varepsilon\Big)\Big)\mathcal{L}w+
\varepsilon b\Big(t,r,w,r\frac{\partial w}{\partial r},\varepsilon\Big)= \nonumber \\
&&=\varepsilon c(t,r,\varepsilon),
\end{eqnarray}
where
\begin{eqnarray*}
a(t,r,w,\Omega,\varepsilon)&=&\varepsilon^{-1}G_I(
\varepsilon y_1+\varepsilon w, \varepsilon v_1+\varepsilon \Omega), \\
b(t,r,w,\Omega,\varepsilon)&=&-(F_I+F_{II})+(F_I+F_{II})\Big|_{w=\Omega=0} \\
c(t,r,\varepsilon)&=&(F_I+F_{II})\Big|_{w=\Omega=0}.
\end{eqnarray*}
Here $v_1$ stands for $r\partial y_1/\partial r$ and
\begin{eqnarray*}
F_I&:=&
-\varepsilon^{-1}G_I(
\varepsilon y_1+\varepsilon w, \varepsilon v_1+\varepsilon\Omega)\mathcal{L}y_1, \\
F_{II}&:=&-\varepsilon^{-2}G_{II}
(r, \varepsilon y_1+\varepsilon w, \varepsilon v_1+\varepsilon\Omega).
\end{eqnarray*}\\

It follows from Proposition 4 that $a,b,c$ are smooth functions of 
$t, x,
(1-x)^{N/2}, w$ and 
$\partial w/\partial x$. Here and hereafter {\bf $x$ denotes the variable
defined by (9), which is equivalently used instead of $r$}.\\

Then the main result of this study can be stated as follows:

\newtheorem{Theorem}{Theorem}
\begin{Theorem}
Assume that $6/5<\gamma \leq 2 (\Leftrightarrow 
4\leq N<12)$ and that $\displaystyle \frac{\gamma}{\gamma-1} (=\frac{N}{2})$
is an integer, that is, $\gamma$ is either $2,3/2,4/3$ or $5/4$.
 Then for any given $T>0$
there is a sufficiently small positive $\varepsilon_0=\varepsilon_0(T)$ such
that, for $|\varepsilon|\leq \varepsilon_0$, there
is a solution $w \in C^{\infty}([0,T]\times[0,R])$ of (13) such that
$$\sup_{j+k\leq n}\Big\|\Big(\frac{\partial}{\partial t}\Big)^j \Big(\frac{\partial}{\partial r}\Big)^k
w
\Big\|_{L^{\infty}([0,T]\times[0,R])} \leq C_n\varepsilon, $$
or a solution $y \in C^{\infty}([0,T]\times[0,R])$ of
(6) or (12) of the form
$$y(t,r)=\varepsilon y_1(t,r)+O(\varepsilon^2), $$
or a motion which can be expressed by the Lagrangian coordinates as
$$r(t,m)=\bar{r}(m)(1+
\varepsilon y_1(t, \bar{r}(m))+ O(\varepsilon^2)) $$
for $0\leq t\leq T, 0\leq m\leq M$.
\end{Theorem}

Our task is to find the inverse image $\mathfrak{P}^{-1}(\varepsilon c)$ of the 
nonlinear mapping $\mathfrak{P}$ defined by
\begin{equation}
\mathfrak{P}(w):=\frac{\partial^2w}{\partial t^2}+(1+\varepsilon a)\mathcal{L}w
+\varepsilon b.
\end{equation}
Note $\mathfrak{P}(0)=0$. It requires a property of the Fr\'{e}chet derivative
of $\mathfrak{P}$:
\begin{equation}D\mathfrak{P}(w)h=
h_{tt}+(1+\varepsilon a_1)\mathcal{L}h+
\varepsilon a_{20}h+
\varepsilon a_{21}rh_r,
\end{equation}
where
\begin{eqnarray*}
a_1(t,r)&=&a\Big(t,r,w,r\frac{\partial w}{\partial r}, \varepsilon\Big), \\
a_{20}(t,r)&=&\frac{\partial a}{\partial w}\mathcal{L}w+\frac{\partial b}{\partial w}, \\
a_{21}(t,r)&=&\frac{\partial a}{\partial\Omega}\mathcal{L}w+\frac{\partial b}{\partial\Omega}.
\end{eqnarray*}
Here $\Omega$ is the dummy of $r\partial w/\partial r$.
We shall use the following observation:

\begin{Proposition}
We have
$$a_{21}=\frac{\gamma P}{\rho}(1+y)^{-2\gamma+2}(1+y+v)^{-\gamma-2}
\Big((\gamma+1)\frac{\partial^2Y}{\partial r^2}+
\frac{4\gamma}{r}\frac{\partial Y}{\partial r}+
\frac{2\varepsilon(\gamma-1)}{1+y}\Big(\frac{\partial Y}{\partial r}\Big)^2\Big),$$
where
$$Y=y_1+w,\qquad y=\varepsilon Y, \qquad v=r\frac{\partial y}{\partial r}=
\varepsilon r\frac{\partial Y}{\partial r}.$$
\end{Proposition}

{\bf Proof}\  Since
\begin{eqnarray*}
\frac{\partial a}{\partial \Omega}&=&\frac{\partial G_I}{\partial v}=
\frac{(1+y)^2}{\gamma}\partial_v^2G_2, \\
\frac{\partial b}{\partial \Omega}&=&\frac{\partial G_I}{\partial v}\mathcal{L}y_1
+\varepsilon^{-1}\frac{\partial G_{II}}{\partial v},
\end{eqnarray*}
we have
$$\varepsilon a_{21}=-(\partial_vG_I)\frac{\gamma P}{\rho }
\Big(\frac{\partial^2y}{\partial r^2}+\frac{4}{r
}\frac{\partial y}{\partial r}\Big)+
\frac{P}{\rho r^2}\partial_vG_{II0}+\frac{1}{\rho r}\frac{dP}{dr}[U],$$
where
$$[U]=-(\partial_vG_I)((3\gamma-4)y+\gamma v)+\partial_vG_{II1}.$$
Since
$$\partial_vG_I=\frac{(1+y)2}{\gamma}\partial_v^2G_2,
\quad
\partial_vG_{II1}=\frac{(1+y)^2}{\gamma}\partial_v^2G_2
((-4+3\gamma)y+\gamma v),
$$
we have $[U]=0$. Using
\begin{eqnarray*}
\partial_v^2G_2&=&-\gamma(\gamma+1)(1+y)^{-2\gamma}(1+y+v)^{-\gamma-2},\\
\partial_vG_{II0}&=&-2\gamma
(1+y)^{-2\gamma+1}(1+y+v)^{-\gamma-2}\cdot
(2(1+y)+(-\gamma+1)v)v,
\end{eqnarray*}
we get the result. \hfill$\blacksquare$\\

\textbullet\  Hereafter we use the variable $x$ defined by
(9) instead of $r=\bar{r}$. \\
 
We note that
$$\frac{\gamma P}{\rho} =\rho_1^{2(\gamma-1)}(1-x)(1+[1-x,(1-x)^{N/2}]_1).$$
Hence the function $\hat{a}_{21}$ defined by
$$\hat{a}_{21}:=
\frac{r}{x(1-x)}\frac{dx}{dr}a_{21}=\frac{r}{\sqrt{x(1-x)}}\rho^{\frac{-\gamma+1}{2}}a_{21}$$
is smooth in $t, x, (1-x)^{N/2}, w, \partial w/\partial x,
\partial^2w/\partial x^2$ including $x=0,1$. Therefore

\begin{Proposition}
The derivative $D\mathfrak{P}$ can be written as
\begin{equation}
D\mathfrak{P}(w)h=
\frac{\partial^2h}{\partial t^2}+
(1+\varepsilon a_1)\mathcal{L}h+
\varepsilon \hat{a}_{21}x(1-x)\frac{\partial h}{\partial x}+
\varepsilon a_{20}h,
\end{equation}
where $a_1, \hat{a}_{21}, a_{20}$ are smooth functions
of $t, x, (1-x)^{N/2}, w, \partial w/\partial x$ and
$\partial^2w/\partial x^2$.
\end{Proposition}

\section{Proof of the main result}

Hereafter we assume that $N/2$ is an integer so that $(1-x)^{N/2}$
is analytic at $x=1$.\\

We are going to apply the Nash-Moser theorem formulated by R. Hamilton (
\cite[p.171, III.1.1.1.]{Hamilton}) as \cite{FE708}, that is: \\

{\bf Nash-Moser(-Hamilton) Theorem}  {\it Let $\mathfrak{E}_0$ and
$\mathfrak{E}$ be tame spaces, $U$ an open subset of $\mathfrak{E}_0$ and $\mathfrak{P}:U \rightarrow \mathfrak{E}$ a smooth tame map. Suppose that the equation for the derivative
$D\mathfrak{P}(w)h=g$ has a unique solution $h=V\mathfrak{P}(w)g$
in $\mathfrak{E}_0$ for all $w$ in $U$ and all $g$ in $\mathfrak{E}$,
and $V\mathfrak{P}:U\times \mathfrak{E}\rightarrow 
\mathfrak{E}_0$ is a smooth tame map. Then $\mathfrak{P}$ is locally invertible.}\\

For the definitions of `tame spaces' and `tame maps', see
\cite{Hamilton} or \cite{FE708}. We shall use the discussions of
\cite{FE708} without repeating the details.\\

We consider the spaces of functions of $t$ and $x$:
\begin{align*}
\mathfrak{E}&:=C^{\infty}([0,T]\times [0,1]) \\
\mathfrak{E}_0&:=\{ w\in \mathfrak{E} \quad |\quad w=\frac{\partial w}{\partial t}=0 \quad\mbox{at}\ t=0\}.
\end{align*}
Let $U$ be the set of all functions $w$ in $\mathfrak{E}_0$ such
that $|w|+|\partial w/\partial x|<1$. Then, for $w \in U$, $y=\varepsilon y_1+\varepsilon w$ and its derivative re $r$ are small, provided that $|\varepsilon|\leq \varepsilon_1$.
Then we can consider the mapping
$$\mathfrak{P}: w\mapsto \partial_t^2w+(1+\varepsilon a)\mathcal{L}w+\varepsilon b$$
maps $U$ into $\mathfrak{E}$, since the coefficients $a, b$ are smooth functions of
$t,x,w,\partial w/\partial x$ and the coefficients
$L_0, L_1$ of $\mathcal{L}$ are analytic on $0\leq x\leq 1$.

The inverse image $\mathfrak{P}^{-1}(\varepsilon c)$ is a desired smooth solution of (13).\\

We should introduce gradings of norms on $\mathfrak{E}$ so that
$\mathfrak{E}, \mathfrak{E}_0$ become  tame spaces in the Hamilton's sense.
To do so,
we use a cut off function
$\omega \in C^{\infty}([0,1])$ such that $\omega(x)=1$ for $0\leq x\leq 1/3$,
$0<\omega(x)<1$ for $1/3<x<2/3$ and $\omega(x)=0$ for $2/3\leq x\leq 1$.
For a function $y$ of $0\leq x\leq 1$, we shall denote
\begin{equation}
y^{[0]}(x)=\omega(x)y(x) , \qquad
y^{[1]}(x)=(1-\omega(x))y(x).
\end{equation} 

We consider the tame spaces
\begin{eqnarray*}
&&\mathfrak{E}_{[0]}=\{y\in C^{\infty}([0,T]\times[0,1])|y=0\ \mbox{for}\  5/6\leq x \leq 1\},\\
&&\mathfrak{E}_{[1]}=\{y\in C^{\infty}([0,T]\times[0,1])|y=0\ \mbox{for}\  
0\leq x \leq 1/6\},
\end{eqnarray*}
endowed with the equivalent gradings of norms
$(\|\cdot\|_{[\mu]n}^{(\infty)})_n, (\|\cdot\|_{[\mu]n}^{(2)})_n, \mu=0,1,$ 
by the same manner as in \cite{FE708},
that is,
denoting
$$\triangle_{[0]}=x\frac{d^2}{dx^2}+\frac{5}{2}\frac{d}{dx}, \qquad
\triangle_{[1]}=z\frac{d^2}{dz^2}+\frac{N}{2}\frac{d}{dz}, \quad (z=1-x),$$
we put
\begin{eqnarray*}
&&\|y\|_{[\mu]n}^{(\infty)}=\sup_{j+k\leq n}\Big\|\Big(-\frac{\partial^2}{\partial t^2}\Big)^j(-\triangle_{[\mu]})^ky\Big\|_{L^{\infty}},\\
&&\|y\|_{[\mu]n}^{(2)}=
\sum_{j+k\leq n}\Big(\int_0^T
\Big\|\Big(-\frac{\partial^2}{\partial t^2}\Big)^j(-\triangle_{[\mu]})^ky\Big\|_{[\mu]}^2dt\Big)^{1/2},
\end{eqnarray*}
where
\begin{eqnarray*}
&&\|y\|_{[0]}=\Big(\int_0^{1}y^2x^{3/2}dx\Big)^{1/2}, \\
&&\|y\|_{[1]}=\Big(\int_{0}^{1}y^2(1-x)^{N/2-1}dx\Big)^{1/2}.
\end{eqnarray*}
On the other hand, on $\mathfrak{E}$ we introduce the gradings of norms
$(\|\cdot\|_n^{(\infty)})_n$ and $(\|\cdot\|_{n}^{(2)})_n$ by
\begin{eqnarray*}
&&\|y\|_n^{(\infty)}:=
\sup_{j+k\leq n, \mu=0,1}\Big\|\Big(-\frac{\partial^2}{\partial t^2}\Big)^j(-\triangle_{[\mu]})^ky^{[\mu]}\Big\|_{L^{\infty}},\\
&&\|y\|_n^{(2)}:=
\Big(\sum_{j+k\leq n, \mu=0,1}\int_0^T\Big\|\Big(-\frac{\partial^2}{\partial t^2}\Big)^j
(-\triangle_{[\mu]})^ky^{[\mu]}\Big\|_{[\mu]}^2dt
\Big)^{1/2}.
\end{eqnarray*}

Then it is easy to see that $\mathfrak{E}$ is a tame space as a tame direct summand of the cartesian
product $\mathfrak{E}_{[0]}\times\mathfrak{E}_{[1]}$, which is a tame space. (See \cite[
p.136, 1.3.3. and 1.3.4.]{Hamilton})
In fact we consider the linear mappings
$L:\mathfrak{E}\rightarrow \mathfrak{E}_{[0]}\times\mathfrak{E}_{[1]}:
h\mapsto (h^{[0]},h^{[1]})$ and
$M:\mathfrak{E}_{[0]}\times\mathfrak{E}_{[1]}:(h_0,h_1)\mapsto h_0+h_1$. It is clear that 
$L$ is tame and $ML=\mbox{Id}_{\mathfrak{E}}$. To verify that $M$ is tame,
we use the following 

\begin{Proposition} If the support of $y(x)$ is included in
$[1/6,5/6]$, then
$$\|\triangle_{[\mu]}^my\|_{L^{\infty}}\leq C
\sum_{0\leq k\leq m}\|\triangle_{[1-\mu]}^ky\|_{L^{\infty}}.$$
\end{Proposition}

A proof can be found in Appendix 2. Now if $h_{\mu}\in \mathfrak{E}_{[\mu]}$,
then $h=M(h_0,h_1)=h_0+h_1$, and
$$h^{[0]}=(h_0+h_1)^{[0]}=\omega h_0+\omega h_1.$$
Then by \cite[Proposition 4]{FE708} we have
$$\|\triangle_{[0]}^mh^{[0]}\|_{L^{\infty}}\leq
C\sum_{k\leq m}\|\triangle_{[0]}^kh_0\|_{L^{\infty}}+\|\triangle_{[0]}^m(\omega h_1)\|_{L^{\infty}}.$$
Proposition 7 can be applied, since
$\mbox{supp}[\omega h_1]\subset [1/6,2/3]$, so that
\begin{eqnarray*}
\mbox{the second term}&\leq& C\sum_{k\leq m}\|\triangle_{[1]}^k(\omega h_1)\|_{L^{\infty}} \\
&\leq& C'\sum_{k\leq m}\|\triangle_{[1]}^kh_{1}\|_{L^{\infty}}.
\end{eqnarray*}
Therefore we have
$$\|\triangle_{[0]}^mh^{[0]}\|_{L^{\infty}}\leq
C\sum_{k\leq m}(\|\triangle_{[0]}^kh_0\|_{L^{\infty}}+
\|\triangle_{[1]}^k
h_1\|_{L^{\infty}}).$$
The same argument gives the estimate of
$\|\triangle_{[1]}^mh^{[1]}\|_{L^{\infty}}$. This implies the tameness of $M$.
Therefore $\mathfrak{E}$ is tame with respect to the grading
$(\|\cdot\|_n^{(\infty)})_n$.\\

By the discussion of \cite{FE708} it is clear that the mapping
$\mathfrak{P}$ is tame. In fact we have
$$\|\mathfrak{P}(w)\|_n^{(\infty)}\leq C
\|w\|_{n+1}^{(\infty)}.$$\\

Therefore we can concentrate ourselves to the 
analysis of the linear equation
\begin{equation}
D\mathfrak{P}(w)h=g
\end{equation}
when $w$ is chosen from $U$ and $g$ is given in $\mathfrak{E}$. By Proposition 3
and 6 we can write 
\begin{equation}
D\mathfrak{P}(w)h=\frac{\partial^2 h}{\partial t^2}-b_2\Lambda h+
b_1x(1-x)\frac{\partial h}{\partial x}+b_0h,
\end{equation}
where
\begin{equation}
\Lambda=x(1-x)\frac{\partial^2}{\partial x^2}+\Big(\frac{5}{2}(1-x)-\frac{N}{2}x\Big)\frac{\partial}{\partial x}
\end{equation}
and
$$b_2=1+\varepsilon a_1, \qquad b_1=
(1+\varepsilon a_1)\frac{L_1}{x(1-x)}+\varepsilon\hat{a}_{21}, $$
$$b_0=(1+\varepsilon a_1)L_0+\varepsilon a_{20}$$
are smooth functions of $t, x, w, Dw, D^2w$, where $D=\partial/\partial x$.\\

In order to establish the existence and uniqueness of the solution of (18),
we introduce the following spaces of
functions of $0\leq x\leq 1$:
\begin{eqnarray*}
\mathfrak{X}=\mathfrak{X}^0&:=&\{y|\|y\|_{\mathfrak{X}}:=
\Big(\int_0^1y^2x^{3/2}(1-x)^{N/2-1}dx
\Big)^{1/2}<\infty\}, \\
\mathfrak{X}^1&:=&\{y\in\mathfrak{X}| \dot{D}y:=\sqrt{x(1-x)}\frac{dy}{dx}\in\mathfrak{X}\}, \\
\mathfrak{X}^2&:=&\{y\in\mathfrak{X}^1|-\Lambda y\in \mathfrak{X}\}.
\end{eqnarray*}

Then we have
\begin{Proposition}
Let $a$ be a function in $C^1[0,1]$. If $y \in \mathfrak{X}^2$ and $v \in \mathfrak{X}^1$, then
$$(-a\Lambda y|v)_{\mathfrak{X}}=(a\dot{D}y|\dot{D}v)_{\mathfrak{X}}+((Da)\check{D}y|v)_{\mathfrak{X}}, $$
where $\check{D}=\displaystyle x(1-x)\frac{d}{dx}$. Here, of course,
$$(u|v)_{\mathfrak{X}}=\int_0^1u v x^{3/2}(1-x)^{N/2-1}dx.$$
\end{Proposition}
{\bf Proof}  If $v\in \mathfrak{X}^1$, then 
$$v(x)=v\Big(\frac{1}{2}\Big)+\int_{\frac{1}{2}}^x
\frac{\dot{D}v(x')}{\sqrt{x'(1-x')}}dx'$$
implies
$$|v(x)|\leq C x^{-3/4}(1-x)^{-N/4+1/2},$$
and
if $y \in \mathfrak{X}^2$, then
\begin{eqnarray*}
x^{5/2}(1-x)^{N/2}\frac{dy}{dx}&=&x^{5/2}(1-x)^{N/2}\frac{dy}{dx}\Big|_{x=1/2} +\\
&&-\int_{1/2}^x\Lambda y(x')x'^{3/2}(1-x')^{N/2-1}dx'
\end{eqnarray*}
implies
$$\Big|x^{5/2}(1-x)^{N/2}\frac{dy}{dx}\Big| \leq C x^{5/4}(1-x)^{N/4}.$$
( Note that the finite constant
$$x^{5/2}(1-x)^{N/2}\frac{dy}{dx}\Big|_{x=1/2}
+\int_0^{1/2}\Lambda y(x')x'^{3/2}(1-x')^{N/2-1}dx'$$
should vanish in order to $\dot{D}y\in\mathfrak{X}$, and so on.)
Therefore the boundary terms in the integration by parts vanish as $x\rightarrow 0, 1$ and we get 
the desired equality. $\blacksquare$\\

Using Proposition 8, we can prove the following energy estimate in the
same manner as \cite[Lemma 3]{FE708}:

\begin{Proposition}
Let $g \in C([0,T], \mathfrak{X})$. If $h\in \bigcap_{k=0,1,2} C^{2-k}([0,T],\mathfrak{X}^k)$
satisfies (18), then we have, for $0\leq t\leq T$,
$$\|\partial_th\|_{\mathfrak{X}}+
\|h\|_{\mathfrak{X}^1}\leq C
(\|\partial_th|_{t=0}\|_{\mathfrak{X}}+
\|h|_{t=0}\|_{\mathfrak{X}^1}+
\int_0^t\|g(t')\|_{\mathfrak{X}}dt').$$
Here
$$\|h\|_{\mathfrak{X}^1}^2=\|h\|_{\mathfrak{X}}^2+\|\dot{D}h\|_{\mathfrak{X}}^2,$$
and the constant $C$ depends only upon
$N, T, \|\partial_tb_2\|_{L^{\infty}},  \|Db_2\|_{L^{\infty}},  \|b_1\|_{L^{\infty}},
\|b_0\|_{L^{\infty}},$ provided that $|1-b_2|\leq 1/2$.
\end{Proposition}

We are considering the initial boundary value problem (IBP):
$$\frac{\partial^2 h}{\partial t^2}+\mathcal{A}h=g(t,x),\qquad h(t,\cdot)\in\mathfrak{X}^1, $$
$$h=\frac{\partial h}{\partial t}=0 \quad \mbox{at}\quad t=0.$$
Here
$$\mathcal{A}=-b_2\Lambda +b_1\check{D}+b_0, \qquad \check{D}=x(1-x)\frac{d}{dx}.$$
Note that ``$h(t,\cdot)\in\mathfrak{X}^1$" is a Dirichlet
boundary condition in some sense. In fact it can be shown that $C_0^{\infty}(0,1)$ is dense in
$\mathfrak{X}^1$. 

Anyway, applying Kato's theory developed in \cite{Kato1970}, we have

\begin{Proposition}
If $g \in C([0,T], \mathfrak{X}^1)\cup C^1([0,T],\mathfrak{X})$, then
there exists a unique solution h of (IBP) in
$\bigcap_{k=0,1,2}C^{2-k}([0,T],\mathfrak{X}^k)$
\end{Proposition}

{\bf Proof}\  We write (IBP) as
$$
\frac{d}{dt}\begin{pmatrix}h \\
\dot{h}\end{pmatrix}+\begin{pmatrix}0 & -1 \\
\mathcal{A} & 0\end{pmatrix}
\begin{pmatrix} h \\
\dot{h}\end{pmatrix}
=\begin{pmatrix}0 \\
g\end{pmatrix}.
$$
Applying the semi-group theory in the space
$\mathfrak{H}=\mathfrak{X}^1\times\mathfrak{X}$ to the family of 
generators
$$D(\mathfrak{A}(t))=\mathfrak{X}^2\times\mathfrak{X}^1, $$
$$\mathfrak{A}(t)=\begin{pmatrix}
0 & -1 \\
\mathcal{A} & 0\end{pmatrix},$$
we get the result. The proof is same as in the Appendix C of
\cite{FE708}. Note that
$$(\mathcal{A}y|v)_{\mathfrak{X}}=(b_2\dot{D}y|\dot{D}v)_{\mathfrak{X}}+(
((b_1+Db_2)\check{D}+b_0)y|v)_{\mathfrak{X}}$$
for $y \in \mathfrak{X}^2$ and $v\in\mathfrak{X}^1$ thanks to Proposition 8.
$\blacksquare$\\

We are going to prove the smoothness of the solution and to get its tame estimates. In order to
do it, we use the cut off function $\omega$ to separate the
singularities at $x=0$ and $x=1$, since, although the singularities 
are of the same type, the calculus structure
of $\Lambda^m, m\in \mathbb{N}$, is little bit complicated.

The equation $\displaystyle \frac{\partial^2 h}{\partial t^2}+\mathcal{A}h=g$ is split into the following
simultaneous system of equations:
\begin{eqnarray}
\Big(\frac{\partial^2}{\partial t^2}+\mathcal{A}_{[0]}\Big)h^{[0]}&=&g^{[0]}-(c_1\check{D}+c_0)h^{[1]} \nonumber \\
\Big(\frac{\partial^2}{\partial t^2}+\mathcal{A}_{[1]}\Big)h^{[1]}&=&g^{[1]}+(c_1\check{D}+c_0)h^{[0]},
\end{eqnarray}
where
\begin{eqnarray*}
&&c_1=(2b_2-b_1)D\omega,\qquad c_0=b_2(\Lambda \omega), \\
&&\mathcal{A}_{[0]}=-b_2\Lambda +(b_1+c_1)\check{D}+b_0+c_0, \\
&&\mathcal{A}_{[1]}=-b_2\Lambda +(b_1-c_1)\check{D}+b_0-c_0.
\end{eqnarray*}

We can rewrite them as:
\begin{eqnarray*}
&&\mathcal{A}_{[0]}=-b_{[0]2}\triangle_{[0]}+b_{[0]1}x\frac{d}{dx}+b_{[0]0}, \\
&&\mathcal{A}_{[1]}=-b_{[1]2}\triangle_{[1]}+b_{[1]1}z\frac{d}{dz}+b_{[1]0}, \qquad (z=1-x),
\end{eqnarray*}
where
\begin{eqnarray*}
&&b_{[0]2}=b_2(1-x),\qquad b_{[1]2}=b_2x, \\
&&b_{[0]1}=\frac{N}{2}b_2+(b_1+c_1)(1-x), \qquad b_{[1]1}=\frac{5}{2}b_2-(b_1-c_1)x, \\
&&b_{[0]0}=b_0+c_0,\qquad b_{[1]0}=b_0-c_0.
\end{eqnarray*}
We may assume that $|b_{[\mu]2}-1|\leq\kappa$ on $x\in I_{[\mu]}, \mu=0,1$, with a constant $\kappa$
such that
 $2/3<\kappa<1$, e.g., $\kappa=5/6$. Here $I_{[0]}=[0,2/3], I_{[1]}=[1/3, 1]$.\\

We note that the regularity of the solution $h$ established by Proposition 10 
can be reduced to that of $h^{[0]}, h^{[1]}$. In fact, if
we know $h^{[0]}\in C^{\infty}([0,T]\times[0,2/3])$, then $h(t,x)=h^{[0]}(t,x)/\omega(x)$
is smooth on $0\leq x< 2/3$, and the smoothness of $h^{[1]}$ implies that of $h(t,x)=h^{[1]}/(1-\omega(x))$ 
on $1/3<x\leq 1$.

But the regularity of the solution of the simultaneous system (21) can be proved by Kato's theory developed 
in \cite{Kato1976}, as in Appendix C of \cite{FE708}. Namely, we consider in the space
\begin{eqnarray*}
\hat{\mathfrak{H}}&=&\mathfrak{H}_{[0]}\times\mathfrak{H}_{[1]}\times \mathbb{R} \\
&=&\mathfrak{X}_{[0]0}^1\times\mathfrak{X}_{[0]}\times
\mathfrak{X}_{[1]0}^1\times\mathfrak{X}_{[1]}\times\mathbb{R}
\end{eqnarray*}
the family of generators
\begin{eqnarray*}
D(\hat{\mathfrak{A}}(t))&=&\hat{\mathfrak{G}}=\mathfrak{G}_{[0]}\times\mathfrak{G}_{[1]}\times \mathbb{R} \\
&=&\mathfrak{X}_{[0](0)}^2\times\mathfrak{X}_{[0]0}^1\times
\mathfrak{X}_{[1](0)}^2\times\mathfrak{X}_{[1]0}^1\times\mathbb{R}, \\
\hat{\mathfrak{A}}(t)&=&\mathfrak{A}_{[0]}(t)\otimes\mathfrak{A}_{[1]}(t)\otimes 0+ B(t),\\
B(t)&=&\begin{pmatrix}0 & 0 & 0 & 0 & 0 \\
0 & 0 & -(c_1\check{D}+c_0) & 0 & -g^{[0]}\\
0 & 0 & 0 & 0 & 0 \\
c_1\check{D}+c_0 & 0 & 0 & 0 & -g^{[1]}\\
0& 0& 0& 0& 0\end{pmatrix},
\end{eqnarray*}
where
$$\mathfrak{A}_{[\mu]}(t)=\begin{pmatrix}0 & -1 \\
\mathcal{A}_{[\mu]} & 0\end{pmatrix}.
$$

Here we set
\begin{eqnarray*}
\mathfrak{X}_{[0]}&=&\{y|\|y\|_{\mathfrak{X}_{[0]}}:=
\Big(\int_0^{2/3}y(x)^2x^{3/2}dx\Big)^{1/2}<\infty\}, \\
\mathfrak{X}_{[0]}^1&=&\{y\in\mathfrak{X}_{[0]} | \dot{D}_{[0]}y=\sqrt{x}\frac{dy}{dx}\in\mathfrak{X}_{[0]}\}, \\
\mathfrak{X}_{[0]0}^1&=&\{y\in\mathfrak{X}_{[0]}^1 |\  y|_{x=2/3}=0\}, \\
\mathfrak{X}_{[0]}^2&=&\{ y\in\mathfrak{X}_{[0]}^1 |\  \triangle_{[0]}y \in \mathfrak{X}_{[0]}\}, \\
\mathfrak{X}_{[0](0)}^2&=&\mathfrak{X}_{[0]}^2\cap\mathfrak{X}_{[0]0}^1;\\
\mathfrak{X}_{[1]}&=&\{y|\|y\|_{\mathfrak{X}_{[1]}}:=
\Big(\int_{1/3}^{1}y(x)^2(1-x)^{N/2-1}dx\Big)^{1/2}<\infty\}, \\
\mathfrak{X}_{[1]}^1&=&\{y\in\mathfrak{X}_{[1]} | \dot{D}_{[1]}y=-\sqrt{1-x}\frac{dy}{dx}\in\mathfrak{X}_{[1]}\}, \\
\mathfrak{X}_{[1]0}^1&=&\{y\in\mathfrak{X}_{[1]}^1 |\  y|_{x=1/3}=0\}, \\
\mathfrak{X}_{[1]}^2&=&\{ y\in\mathfrak{X}_{[1]}^1 |\  \triangle_{[1]}y \in \mathfrak{X}_{[1]}\}, \\
\mathfrak{X}_{[1](0)}^2&=&\mathfrak{X}_{[1]}^2\cap\mathfrak{X}_{[1]0}^1.
\end{eqnarray*}\\

{\bf Remark }\  1) It may be
difficult to verify that,
given a solution $(h_0,h_1)$ of the system (21)
such that $h_{\mu}\in \bigcap_{k=0,1,2}C^{2-k}([0,T], \mathfrak{X}_{[\mu]}^k)$, the function
$h$ which should be defined by
$$h(t,x)=\begin{cases}
h_0(x) & \qquad (0\leq x\leq 1/3) \\
h_0(x)+h_1(x) & \qquad (1/3< x<2/3) \\
h_1(x) & \qquad (2/3\leq x\leq 1)
\end{cases}
$$
belongs to $C([0,T], \mathfrak{X}^2)$. Therefore we first established
the existence of the solution $h$ by Proposition 10. Then, by the uniqueness, 
we can claim that $h^{[\mu]}=h_{\mu}$, the solutions of (21).

2) We used
$$\|y\|_{[0]}=\Big(\int_0^1y(x)^2x^{3/2}dx\Big)^{1/2}$$
in the definition of the gradings on
$\mathfrak{E}_{[0]}$. But 
$\|y\|_{\mathfrak{X}_{[0]}}=\|y\|_{[0]}$ for $y=h^{[0]}$, since
$\mbox{supp}[h^{[0]}]\subset [0, 2/3]$. So, we can consider
$h(t,\cdot)^{[\mu]}\in \mathfrak{X}_{[\mu](0)}^2$ for the solution $h$
established in Proposition 10.\\

Then $B(t) \in C([0,T],\mathsf{B}(\hat{\mathfrak{H}}))$ is a smooth bounded perturbation from the
stable family $(\mathfrak{A}_{[0]}(t)\otimes\mathfrak{A}_{[1]}(t)\otimes 0)_t$.
Hence $(\hat{\mathfrak{A}}(t))_t$ is stable.\\

In order to consider `smoothness', `ellipticity' and compatibility conditions, we introduce the scales of Hilbert spaces 
$$ \hat{\mathfrak{H}}_j=\mathfrak{X}_{[0](0)}^{j+1}\times
\mathfrak{X}_{[0]}^j\times\mathfrak{X}_{[1](0)}^{j+1}\times
\mathfrak{X}_{[1]}^j\times\mathbb{R}, $$
$$\hat{\mathfrak{G}}_j=\hat{\mathfrak{G}}\cap\hat{\mathfrak{H}}_j=
\mathfrak{X}_{[0](0)}^{j+1}\times
\mathfrak{X}_{[0](0)}^j\times\mathfrak{X}_{[1](0)}^{j+1}\times
\mathfrak{X}_{[1](0)}^j\times\mathbb{R},$$
as in Appendix C of \cite{FE708}, where
\begin{align*}
\mathfrak{X}_{[\mu]}^{2m+1}&=
\{ y\in\mathfrak{X}_{[\mu]}^{2m}\quad|\quad \dot{D}_{[\mu]}\triangle_{[\mu]}^m
y\in\mathfrak{X}_{[\mu]}\}, \\
\mathfrak{X}_{[\mu]}^{2m+2}&=
\{y\in\mathfrak{X}_{[\mu]}^{2m+1}\quad |\quad \triangle_{[\mu]}^{m+1}y
\in\mathfrak{X}_{[\mu]}\},\\
\mathfrak{X}_{[\mu](0)}^j&=\mathfrak{X}_{[\mu]}^j\cap
\mathfrak{X}_{[\mu]0}^1.
\end{align*}
The definition of
$\|\cdot\|_{\mathfrak{X}_{[\mu]}^j}$ follows that of
$\|\cdot\|_j$ in \cite{FE708}, that is:
$$\|y\|_{\mathfrak{X}_{[\mu]}^j}=\Big(
\sum_{\ell\leq j} (\langle y\rangle_{[\mu]\ell})^2\Big)^{1/2}, $$
$$\langle y\rangle_{[\mu]\ell}=
\begin{cases}
\|\triangle_{[\mu]}^my
\|_{\mathfrak{X}_{[\mu]}}
 \quad\mbox{as}\quad \ell=2m, \\
\|\dot{D}_{[\mu]}\triangle_{[\mu]}^my\|_{\mathfrak{X}_{[\mu]}} \quad\mbox{as}
\quad
\ell=2m+1.
\end{cases}
$$

In order check the `smoothness', we note that
$c_1=c_0=0$ for $0\leq x\leq 1/3$ or $2/3\leq x\leq 1$. This implies that
\begin{eqnarray*}
&&\|(c_1\check{D}+c_0)y^{[1]}\|_{\mathfrak{X}_{[0]}^j}\leq C\|(c_1\check{D}+c_0)y^{[1]}\|_{\mathfrak{X}_{[1]}^j}
\leq C'\|y^{[1]}\|_{\mathfrak{X}_{[1]}^{j+1}}, \\
&&\|(c_1\check{D}+c_0)y^{[0]}\|_{\mathfrak{X}_{[1]}^j}\leq C |(c_1\check{D}+c_0)y^{[0]}\|_{\mathfrak{X}_{[0]}^j}\leq
C'\|y^{[0]}\|_{\mathfrak{X}_{[0]}^{j+1}}.
\end{eqnarray*}
(See \cite[Proposition 6]{FE708}.) 
Here we have used the following

\begin{Proposition}
If the support of $y \in C^{\infty}(0,1)$ is included in $[1/3,2/3]$,
then 
$$\|y\|_{\mathfrak{X}_{[\mu]}^j}\leq C \|y\|_{\mathfrak{X}_{[1-\mu]}^j},$$
where $\mu=0,1$.
\end{Proposition}
A proof can be found in Appendix 2.

Then, using this observation, we can reduce the `ellipticity' of
$\hat{\mathfrak{A}}(t)$ to that of
$\mathcal{A}_{[\mu]}(t), \mu=0,1$.\\

The compatibility conditions are guaranteed as follows.

We are considering the Cauchy problem
$$\frac{du}{dt}+\hat{\mathfrak{A}}(t)u=0,\qquad u|_{t=0}=\phi_0, $$
where
\begin{align*}
\hat{\mathfrak{A}}(t)&=\begin{pmatrix}0 & -1 & 0 & 0 & 0 \\
\mathcal{A}_{[0]}& 0 &-\mathcal{C} & 0 &-g^{[0]} \\
0 & 0 & 0 & -1 & 0 \\
\mathcal{C} & 0 & \mathcal{A}_{[1]} & 0 & -g^{[1]} \\
0 & 0 & 0 & 0 & 0 \end{pmatrix}, \\
\mathcal{C}& :=c_1\check{D}+c_0, \\
\phi_0&=\begin{pmatrix}
0\\
0\\
0\\
0\\
1\end{pmatrix}.
\end{align*}
As in \cite[Section 2]{Kato1970}, we consider
\begin{align*}
S^0&=I, \\
S^{j+1}\phi&=-\sum_{k=0}^j\binom{j}{k}\Big(\frac{d}{dt}\Big)^{j-k}
\hat{\mathfrak{A}}(0)S^k\phi, \\
D_0&=\hat{\mathfrak{H}}=\mathfrak{X}_{[0]0}^1\times
\mathfrak{X}_{[0]}\times\mathfrak{X}_{[1]0}^1\times
\mathfrak{X}_{[1]}\times\mathbb{R}, \\
D_{j+1}&=\{\phi\in D_j| S^k\phi
\in \hat{\mathfrak{G}}_{j+1-k}, \quad 0\leq k\leq j\}.
\end{align*}

We should show that $\phi_0\in D_n$ for any $n$. But
$g^{[0]}, g^{[1]}$ can be considered as functions in 
$C^{\infty}([0,T]\times[0,1])$ such that, for all positive integer $\ell$,
$\partial_t^{\ell}g^{[0]}(0,x)=0$ for $2/3\leq x\leq 1$ and
$\partial_t^{\ell}g^{[1]}(0,x)=0$ for $0\leq x\leq 1/3$. We denote
$$\phi_k:=S^k\phi_0=\begin{pmatrix}\phi_{[0]0}^k\\
\phi_{[0]1}^k\\
\phi_{[1]0}^k\\
\phi_{[1]1}^k\\
0\end{pmatrix}.$$
Then it is easy to verify by induction that, for $k\geq 1$,
the extension $\tilde{\phi}_k
=(\tilde{\phi}_{[0]0}^k,\tilde{\phi}_{[0]1}^k, \tilde{\phi}_{[1]0}^k,
\tilde{\phi}_{[1]1}^k, 0)^T$ of $\phi_k$ defined by
\begin{align*}
\tilde{\phi}_{[0]0}^k(x)&=\begin{cases}
\phi_{[0]0}^k(x) &\quad (0\leq x\leq 2/3) \\
0 &\quad (2/3<x\leq 1)
\end{cases}\\
\tilde{\phi}_{[0]1}^k(x)&=\begin{cases}
\phi_{[0]1}^k(x) &\quad (0\leq x\leq 2/3) \\
0 &\quad (2/3<x\leq 1)
\end{cases}\\
\tilde{\phi}_{[1]0}^k(x)&=\begin{cases}
0 &\quad (0\leq x<1/3)\\
\phi_{[1]0}^k(x) &\quad (1/3\leq x\leq 1) 
\end{cases}\\
\tilde{\phi}_{[1]1}^k(x)&=\begin{cases}
0 &\quad(0\leq x<1/3) \\
\phi_{[1]1}^k(x) &\quad (1/3\leq x\leq 1) 
\end{cases}.
\end{align*}
belongs to $C^{\infty}([0,1];\mathbb{R}^5)$.
In other words, the components of $\phi_k$ satisfy the boundary conditions at $x=1/3$ and $x=2/3$ and 
$\phi_k=S^k\phi_0$ remains in $\hat{\mathfrak{G}}_{k+1}$. It implies that $\phi_0\in D_n$ for all
$n$. 

Summing up, we can claim that $h^{[0]}\in C^{\infty}([0,T]\times[0,2/3])$ and
$h^{[1]}\in C^{\infty}([0,T]\times [1/3,1])$ provided that
$g\in C^{\infty}([0,T]\times[0,1])$.\\

Finally, we must deduce the tame estimate of
$(w,g)\mapsto h$. We are going to show that
$$\|h\|_{n+2}^{\langle T\rangle}\leq C
(1+\|g\|_{n+1}^{\langle T\rangle}+|w|_{n+7}^{\langle T\rangle}).$$
Here
\begin{eqnarray*}
\|y\|_n^{\langle T\rangle}&:=&\Big(
\sum_{j+k\leq n, \mu=0,1}\int_0^T\|\partial_t^jy^{[\mu]}\|_{\mathfrak{X}_{[\mu]}^k}dt\Big)^{1/2}, \\
|y|_n^{\langle T\rangle}&:=&\max_{j+k\leq n,\mu=0,1}
\|\partial_t^j\dot{D}_{[\mu]}^ky^{[\mu]}\|_{L^{\infty}([0,T]\times[0,1])}.
\end{eqnarray*}\\

Let us follow the discussion of \cite[\S 5.4]{FE708}. To do so, 
we should reconsider the discussion about the single equation, say,
we consider a solution $H$ of the boundary value problem
$$\frac{\partial^2H}{\partial t^2}+\mathcal{A}(\vec{b})H=G(t,x),\quad H|_{x=1}=0$$
on $0\leq t\leq T$. Here $\vec{b}$ stands for the vector
$(b_0,b_1,b_2)$. The energy estimate claimed in Proposition 9 should read
$$\|\partial_tH\|+\|H\|_1\leq C(
\|\partial_tH|_{t=0}\|+
\|H|_{t=0}\|_1+
\int_0^t\|G(t')\|dt').$$
Even if we consider the $H=h$ which satisfy the initial condition
$h|_{t=0}=\partial_th|_{t=0}=0$, the higher derivatives
$\partial_t^{n+2}h$ may not vanish at $t=0$. Therefore the estimate
of $\|\partial_t^{n+1}h\|_1$ in the proof of
\cite[Proposition 10]{FE708}
should be replaced by
\begin{align*}
\|\partial_t^{n+1}h\|_1&\leq C(\|\partial_t^{n+2}h|_{t=0}\|+
\|\partial_t^{n+1}h|_{t=0}\|_1+ \\
&+\int_0^t\|\partial_t^{n+1}g\|dt'+\int_0^t\|[\partial_t^{n+1},\mathcal{A}]h\|dt').
\end{align*}
We claim the estimate
\begin{equation}
\|\partial_t^{n+2}h|_{t=0}\|+
\|\partial_t^{n+1}h|_{t=0}\|_1\leq
C(1+W_n(g)+
|\vec{b}|_{n+1}^{\langle 0\rangle}),
\end{equation}
provided that $ W_0(g), |\vec{b}|_4^{\langle 0\rangle}\leq
M_0$. Here
$$W_n(g):=\sum_{j+k\leq n}\|\partial_t^jg|_{t=0}\|_k$$
and
$$|y|_n^{\langle 0\rangle }:=\max_{j+k\leq n}
\|\partial_t^j\dot{D}^ky|_{t=0}\|_{L^{\infty}([0,1])}.$$
To prove (22) it is sufficient to verify the following estimate by induction on $n$:
for all $k\in \mathbb{N}$,
$$\|\partial_t^{n+2}h|_{t=0}\|_k
\leq C(
|\vec{b}|_{n+k+1}^{\langle 0\rangle}W_0(g)+
|\vec{b}|_{k+3}^{\langle 0\rangle}W_{n-2}(g)+W_{n+k}(g)).
$$
Since the proof of the above inequality by induction on $n$ using
the estimate
$$\|\mathcal{A}(\vec{b})y\|_k\leq C(\|y\|_{k+2}+ |\vec{b}|_{k+3}\|y\|)$$
applied to the relation
$$\partial_t^{n+2}h=-\sum_{j=0}^n
\binom{n}{j}\mathcal{A}(\partial_t^{n-j}\vec{b})\partial_t^jh+\partial_t^ng$$
is straightforward, we omit it. 

Moreover we note that the inequality in the statement of
\cite[Lemma 4]{FE708} can be replaced by the stronger one:
$$\|h\|_{n+2}^{\langle t \rangle}\leq
C\Big(1+\int_0^t
\|g\|_{n+1}^{\langle t' \rangle}dt'+
W_n(g)+\|g\|_n^{\langle T \rangle}+|\vec{b}|_{n+3}^{\langle T \rangle}\Big),$$
for $0\leq t\leq T$, where
\begin{eqnarray*}
\|y\|_n^{\langle\tau\rangle}&&:=\Big(
\sum_{j+k\leq n}\int_0^{\tau}\|\partial_t^jy\|_{k}^2dt\Big)^{1/2}, \\
|y|_n^{\langle\tau\rangle}&&:=\max_{j+k\leq n}\|\partial_t^j\dot{D}^ky\|_{L^{\infty}([0,\tau]\times[0,1])}
\end{eqnarray*}
This can be verified easily by following the discussion in \cite[\S 5.4]{FE708}.
Let us omit the detail. \\

Let us go back to the simultaneous system of equations. 
Applying the above discussion on a single equation, 
we have
\begin{eqnarray*}
\|h^{[0]}\|_{[0]n+2}^{\langle t \rangle}\leq C\Big(1+\int_0^t\|h^{[1]}\|_{[1]n+2}^{\langle t' \rangle}dt'+
W_n(g)+
\|g\|_{n+1}^{\langle T \rangle}+|\vec{b}|_{n+3}^{\langle T \rangle}\Big), \\
\|h^{[1]}\|_{[1]n+2}^{\langle t \rangle}\leq C\Big(1+\int_0^t\|h^{[0]}\|_{[0]n+2}^{\langle t' \rangle}dt'
+
W_n(g)+\|g\|_{n+1}^{\langle T \rangle}+|\vec{b}|_{n+3}^{\langle T \rangle}\Big), 
\end{eqnarray*}
for $0\leq t \leq T$, since
$$\|(c_1\check{D}+c_0)h^{[\mu]}\|_{[1-\mu]k}\leq C(1+\|h^{[\mu]}\|_{[\mu]k+1}+
|\vec{b}|_{k+3}^{\langle T\rangle})$$
for $\mu=0,1$. Here $\|\cdot\|_{[\mu]k}$ stands for
$\|\cdot\|_{\mathfrak{X}_{[\mu]}^k}$. 
Applying the Gronwall's lemma to
the quantity
$$U(t):=\|h^{[0]}\|_{[0]n+2}^{\langle t\rangle}+\|h^{[1]}\|_{[1]n+2}^{\langle t\rangle},$$
we get
$$U(t)\leq C(1+W_n(g)+\|g\|_{n+1}^{\langle T\rangle}
+|\vec{b}|_{n+3}^{\langle T\rangle}).$$

This completes the proof, since $W_n(g)\leq C\|g\|_{n+1}^{\langle T\rangle}$ by Sobolev's imbedding.

\section{Cauchy problems}

We have discussed about the justification of
linearized approximations by time-periodic solutions. In this section
we want to give a brief mention on the Cauchy problems associated with
the equation (6) or (12). We consider the problem (CP):
\begin{align*}
&\frac{\partial^2y}{\partial t^2}+\Big(1+G_I\Big(y,r\frac{\partial y}{\partial r}\Big)\Big)\mathcal{L}y+
G_{II}\Big(r,y,r\frac{\partial y}{\partial r}\Big)=0, \\
&y|_{t=0}=\psi_0(r),\qquad \frac{\partial y}{\partial t}\Big|_{t=0}=\psi_1(r),
\end{align*}
where the initial data $\psi_0, \psi_1$ are given functions.
We claim

\begin{Theorem}
Assume that $6/5<\gamma \leq 2 (\Leftrightarrow 
4\leq N<12)$ and that $\displaystyle \frac{\gamma}{\gamma-1} (=\frac{N}{2})$
is an integer, that is, $\gamma$ is either $2,3/2,4/3$ or $5/4$.
Then for any given $T>0$ there exist a sufficiently small positive number $\delta$ and a sufficiently large integer 
$\mathfrak{r}$ such that if
$\psi_0,\psi_1 \in C^{\infty}([0,R])$ satisfy
$$\max_{j\leq2(2\mathfrak{r}+1)}\Big\{\Big\|\Big(\frac{d}{dr}\Big)^j\psi_0\Big\|_{L^{\infty}(0,R)},
\Big\|\Big(\frac{d}{dr}\Big)^j\psi_1\Big\|_{L^{\infty}(0,R)}\Big\}\leq\delta,
$$
then there exists a unique solution $y(t,r)$ of 
(CP) in $C^{\infty}([0,T]\times [0,R])$.
\end{Theorem}

A proof of this theorem can be done as follows.

Let us take the function
$$y_1^*(t,r)=\psi_0(r)+t\psi_1(r),$$
which satisfy the initial conditions. Then we should find a solution $w$ introduced by
$$y=y_1^*+w,$$
which should obey the initial conditions
$$w|_{t=0}=\frac{\partial w}{\partial t}\Big|_{t=0}=0.$$
The equation which $w$ should satisfies is same as (13), in which the
time-periodic function 
$$\varepsilon y_1=\sin(\sqrt{\lambda}t+\theta)\Phi(r)$$
is replaced by
$$y_1^*=\psi_0(r)+t\psi_1(r),$$
and $F_I+F_{II}$ should be replaced by
$$(1+G_I(y_1^*+w,v_1^*
+\Omega))\mathcal{L}y_1^*+G_{II}(r,y_1^*+w, v_1^*+\Omega).$$
Of course we take $\varepsilon =1$. Then the mapping
$\mathfrak{P}(w)$ and the derivative $D\mathfrak{P}(w)h$ are
defined in the same forms as (14) and as (15). 
Proposition 5 holds valid, since the concrete structure of the function $y_1$ or
$y_1^*$ is not used in the proof; It is sufficient that $\varepsilon y_1$ or $y_1^*$ is a
small smooth function. Hence Proposition 6 holds valid, when $\varepsilon y_1$ is replaced by $y_1^*$.

Then the proof of Theorem 1 given in \S 4 can be repeated word for word 
in the present situation.
Note that
$$c=-\Big(1+G_I\Big(y_1^*,r\frac{\partial y_1^*}{\partial r}\Big)\Big)\mathcal{L}y_1^*
-G_{II}\Big(r,y_1^*,r\frac{\partial y_1^*}{\partial r}\Big)$$
and
that $\|c\|_n^{(\infty)}\leq C(\|\psi_0\|_{n+1}^{(\infty)}+\|\psi_1\|_{n+1}^{(\infty)})$,
provided that $0\leq t\leq T$. In fact, if we follow the discussion of
\cite[III.1.]{Hamilton}, we can show that it is enough to take
$\mathfrak{r}$ such that $2\mathfrak{r}>3/2+\max\{5,N\}/4$. 
(But this $\mathfrak{r}$ may not be the best possible.) Anyway
this completes the proof of Theorem 2.\\

{\bf Remark}\  The corresponding initial data in the Eulerian variables are given by
\begin{align*}
\rho|_{t=0}(r)&=\bar{\rho}(\bar{r})\Big((1+\psi_0(\bar{r}))^2
\Big(1+\psi_0(\bar{r})+\bar{r}\frac{d\psi_0(\bar{r})}{d\bar{r}}\Big)\Big)^{-1}, \\
u|_{t=0}(r)&=\bar{r}\psi_1(\bar{r})
\end{align*}
implicitly by $\bar{r}=\bar{r}(m(r))$. Here $m\mapsto \bar{r}(m)$ is the 
inverse function of
$$\bar{r}\mapsto m=m(\bar{r})=4\pi\int_0^{\bar{r}}\bar{\rho}(r)r^2dr$$ and
$r\mapsto m(r)$ is the inverse function of
$m\mapsto r=\bar{r}(m)(1+\psi_0(\bar{r}(m))$.

\section{Concluding remark}

In order that the equilibrium satisfy
that $\bar{\rho}^{\gamma-1}$ is analytic at the free boundary $r=R$ and that 
the eigenfunction $y_1$ is analytic in $r$ at $r=R$, we have assumed 
that $N$ is an even integer. But $\gamma=5/3 (N=5)$ for mono-atomic gas, and
$\gamma=7/5 (N=7)$ for the air. Therefore
it is desired that the result will
be extended to the case when $N$ is an odd integer at least. Moreover
for the case when $N$ is not an integer, we might try quite other approach.
It seems that these are interesting open problems
in view of physical applications. \\

{\bfseries\Large Acknowledgment} The author would like to express
his sincere and deep thanks to the anonymous referee
for his/her careful reading of the original manuscript, which contained many careless typos and rude expressions. If this revised manuscript has turned out to be tolerably readable, the author owes it to the kind suggestions by the referee.\\

{\bfseries\Large Appendix 1 }\\

Let us consider a solution 
$\rho =\rho(r), r_0\leq r<R,$ of the Lane-Emden
equation
$$-\frac{1}{r^2}\frac{d}{dr}\Big(\frac{r^2}{\rho}\frac{dP}{dr}\Big)=
4\pi g_0\rho, \quad P=A\rho^{\gamma}. $$
Let $[r_0,R)$ be a right maximal interval of existence of
$\rho >0$, and we assume that $R<+\infty, d\rho/dr|_{r=r_0}<0$.
Then there is a positive constant $C$ such that
$$\rho=C(R-r)^{\frac{1}{\gamma-1}}\Big(
1+\Big[\frac{R-r}{R}, C'\Big(\frac{R-r}{R}\Big)^{\frac{\gamma}{\gamma-1}}\Big]_1\Big) $$
with
$$C'=KR^{\frac{\gamma}{\gamma-1}}
C^{2-\gamma}, \qquad K=\frac{4\pi g_0(\gamma-1)}{A\gamma}.$$\\

{\bf Proof}\  The variable
$$U:=\rho^{\gamma-1}$$
satisfies
$$\frac{d^2U}{dr^2}+\frac{2}{r}\frac{dU}{dr}+KU^m=0, $$
where
$m=1/(\gamma-1)$.
Then
$$v:=-\frac{r}{U}\frac{dU}{dr}, \qquad w:=Kr^2U^{m-1}$$
satisfies the plane autonomous system
\begin{eqnarray*}
r\frac{dv}{dr}&=&-v+v^2+w \\
r\frac{dw}{dr}&=&w(2-(m-1)v)
\end{eqnarray*}

The interval $[r_0,R)$ is right maximal. We assumed that $v(r_0)>0$.
We claim that there is $r_1 \in [r_0,R)$ such that
$v(r_1)>1$. Otherwise
$0<v\leq 1$ and $|\frac{r}{w}\frac{dw}{dr}|\leq m+1$ for $r_0\leq r <R$.
Then it should be $R=+\infty$, a contradiction to the assumption.
Hence we can assume that $v(r_0)>1$. Then
$r\frac{dv}{dr}\geq v(v-1)$ implies $v\geq 1+\delta$, $dv/dr >0$ and 
$r\frac{dw}{dr}\leq 2w$. So, it should be that $v(r)\rightarrow
+\infty $ as $r\rightarrow R$, since $R<\infty$. We see $w \leq B$.

Now we introduce the variables
$$x_1:=\frac{1}{v}, \qquad x_2:=\frac{w}{v^2}, $$
$$t:=\exp \Big(-\int_{r_0}^r
\frac{v(r')dr'}{r'}\Big). $$
Then $(x_1,x_2) \rightarrow (0,0), t\rightarrow 0$ as
$r \rightarrow R$ and $(x_1(t), x_2(t)), 0<t\leq 1,$ satisfies
\begin{eqnarray*}
t\frac{dx_1}{dt}&=&(1-x_1+x_2)x_1 \\
t\frac{dx_2}{dt}&=&(m+1-4x_1+2x_2)x_2.
\end{eqnarray*}
As well-known, this Briot-Bouquet system can be reduced to
\begin{eqnarray*}
t\frac{dz_1}{dt}&=&z_1 \\
t\frac{dz_2}{dt}&=&(m+1)z_2
\end{eqnarray*}
by a transformation of the form
\begin{eqnarray*}
x_1&=&z_1(1+P_1(z_1,z_2)) \\
x_2&=&z_2(1+P_2(z_1,z_2)).
\end{eqnarray*}
Here
$$P_j(z_1,z_2)=[z_1,z_2]_1$$
for $j=1,2$. Therefore there are positive constants
$C_1, C_2$ such that
\begin{eqnarray*}
x_1&=& C_1t(1+P_1(C_1t, C_2t^{m+1})), \\
x_2&=& C_2t^{m+1}(1+P_2(C_1t, C_2t^{m+1})).
\end{eqnarray*}
Since $dr/r=-x_1dt/t$, we see
\begin{eqnarray*}
\log \frac{R}{r}&=&\frac{R-r}{R}\Big(1+\Big[\frac{R-r}{R}\Big]_1\Big) \\
&=&C_1t(1+[C_1t, C_2t^{m+1}]_1),
\end{eqnarray*}
from which
$$C_1t=\frac{R-r}{R}\Big(
1+\Big[\frac{R-r}{R},
C'\Big(\frac{R-r}{R}\Big)^{m+1}\Big]_1\Big)$$
and
$$x_1=\frac{R-r}{R}\Big(1+\Big[\frac{R-r}{R},
C'\Big(\frac{R-r}{R}\Big)^{m+1}\Big]_1\Big),$$
where $C'=C_2/C_1^{m+1}$.
Integrating
$dU/U=-dr/rx_1$, we have
$$U=C_3
\frac{R-r}{R}\big(
1+\Big[\frac{R-r}{R},
C'\Big(\frac{R-r}{R}\Big)^{m+1}\Big]_1\Big).$$ 
It is easy to see $C'=KR^2C_3^{m-1}$, and we get the required result.
$\blacksquare$\\

{\bfseries\Large Appendix 2}\\

\textbullet\  Let us prove Proposition 7, that is, 
$$\|\triangle_{[0]}^my\|_{L^{\infty}}\leq C
\sum_{k\leq m}\|\triangle_{[1]}^ky\|_{L^{\infty}},$$
provided that $\mbox{supp}[y]\subset [1/6,5/6]$.

Note that
$$\triangle_{[0]}=\alpha\triangle_{[1]}+\beta\check{D}_{[1]},$$
where $\check{D}_{[1]}=zd/dz=-(1-x)d/dx$ and
$$\alpha=\frac{x}{1-x}, \qquad \beta=-\frac{1}{1-x}\Big(\frac{x}{1-x}\frac{N}{2}+\frac{5}{2}
\Big)$$
are smooth function on $(0,1)$. Therefor our task is
to estimate
$$
\|(\alpha\triangle_{[1]}+\beta\check{D}_{[1]})^my\|_{L^{\infty}}.$$
On the other hand, it is easy to verify that there are
$\gamma_{\epsilon k}^{(m)}\in C^{\infty}(0,1)$ such that
$$(\alpha\triangle_{[1]} +\beta \check{D}_{[1]})^m=
\sum_{k\leq m}(\gamma_{1k}^{(m)}\check{D}_{[1]}\triangle_{[1]}^k+
\gamma_{0k}^{(m)}\triangle_{[1]}^k)$$
with $\gamma_{1m}^{(m)}=0$. 
Note that
$$\|\check{D}_{[1]}\triangle_{[1]}^ky\|_{L^{\infty}}\leq
\|D\triangle_{[1]}^k\|_{L^{\infty}}\leq \frac{2}{N}
\|\triangle_{[1]}^{k+1}y\|_{L^{\infty}}.$$
(See \cite[Proposition 3]{FE708}). This completes the proof.\\

\textbullet\  Let us prove Proposition 11, that is, 
$$\|y\|_{\mathfrak{X}_{[0]}^j}\leq C\|y\|_{\mathfrak{X}_{[1]}^j},$$
provided that $\mbox{supp}[y]\subset [1/3,2/3]$.

It is clear that 
$$\|y\|_{\mathfrak{X}_{[0]}}\leq C\|y\|_{\mathfrak{X}_{[1]}},$$
since $x^{3/2}\leq 3^{N/2-1}(1-x)^{N/2-1}$ for $1/3\leq x\leq 2/3$.
Let us estimate
$\|\triangle_{[0]}^my\|_{\mathfrak{X}_{[0]}}$ and
$\|\dot{D}_{[0]}\triangle_{[0]}^my\|_{\mathfrak{X}_{[0]}}$,
where $\dot{D}_{[0]}=\sqrt{x}d/dx$. As in the above discussion we note that
$$\triangle_{[0]}=\alpha\triangle_{[1]}+\beta\check{D}_{[1]},$$
where $\check{D}_{[1]}=zd/dz=-(1-x)d/dx$ and
$$\alpha=\frac{x}{1-x}, \qquad \beta=-\frac{1}{1-x}\Big(\frac{x}{1-x}\frac{N}{2}+\frac{5}{2}
\Big)$$
are smooth function on $(0,1)$. Therefor our task is
to estimate
$$\|\triangle_{[0]}^my\|_{\mathfrak{X}_{[0]}}\leq C\|
\triangle_{[0]}^my\|_{\mathfrak{X}_{[1]}}=
C\|(\alpha\triangle_{[1]}+\beta\check{D}_{[1]})^my\|_{\mathfrak{X}_{[1]}}$$
and
$$\|\dot{D}_{[0]}\triangle_{[0]}^my\|_{\mathfrak{X}_{[0]}}\leq C\|\check{D}_{[1]}
\triangle_{[0]}^my\|_{\mathfrak{X}_{[1]}}=
C\|\check{D}_{[1]}(\alpha\triangle_{[1]}+\beta\check{D}_{[1]})^my\|_{\mathfrak{X}_{[1]}}.$$

On the other hand, it is easy to verify that there are
$\gamma_{\epsilon k}^{(m)}, \gamma_{\epsilon k}^{(m)\sharp}\in C^{\infty}(0,1)$ such that
\begin{eqnarray*}
&&(\alpha\triangle +\beta \check{D})^m=
\sum_{k\leq m}(\gamma_{1k}^{(m)}\check{D}\triangle^k+
\gamma_{0k}^{(m)}\triangle^k), \\
&&\check{D}(\alpha\triangle+\beta\check{D})^m=
\sum_{k\leq m}(\gamma_{1k}^{(m)\sharp}\check{D}\triangle^k
+\gamma_{0k}^{(m)\sharp}\triangle^k)
\end{eqnarray*}
with $\gamma_{1m}^{(m)}=0$. Here $\triangle, \check{D}$ stand for 
$\triangle_{[1]}, \check{D}_{[1]}$.
Hence we have
\begin{eqnarray*}
&&\|\triangle_{[0]}^my\|_{\mathfrak{X}_{[0]}}\leq C\|y\|_{\mathfrak{X}_{[1]}^{2m}},\\
&&\|\dot{D}_{[0]}\triangle_{[0]}^my\|_{\mathfrak{X}_{[0]}}\leq C\|y\|_{\mathfrak{X}_{[1]}^{2m+1}}.
\end{eqnarray*}
This completes the proof.

\end{document}